\documentclass[fleqn]{amsart}
\usepackage{amssymb, graphics, amsmath, epsfig, changebar}
\usepackage{mathtools} 
\usepackage[percent]{overpic}
\usepackage[caption=false]{subfig} 

\theoremstyle{plain}
\newtheorem{theorem}{Theorem}
\newtheorem{corollary}{Corollary}

\newtheorem{lemma}{Lemma}
\newtheorem{proposition}{Proposition}
\newtheorem{remark}{Remark}

\theoremstyle{definition}

\theoremstyle{example}
\newtheorem{example}{\rm{Example}}

\theoremstyle{remark}

\numberwithin{equation}{section}
\newcommand{\N}{\mathbb{N}}
\newcommand{\Z}{\mathbb{Z}}
\newcommand{\R}{\mathbb{R}}
\newcommand{\C}{\mathbb{C}}
\newcommand{\Se}{S^3_\epsilon}
\newcommand{\Sum}{\displaystyle\sum}

\newcommand{\Om}{\Omega}

\graphicspath{{./images/}}

\newcommand{\overimage}[3]{\begin{overpic}{#1}\put(#2){#3}\end{overpic}}
\newcommand{\oversimage}[5]{\begin{overpic}{#1}\put(#2){#3}\put(#4){#5}\end{overpic}}

\begin{document}                   

\title{Knots with Hopf crossing number at most one}
\author{Maciej Mroczkowski}
\address{Institute of Mathematics\\
Faculty of Mathematics, Physics and Informatics\\
University of Gdansk, 80-308 Gdansk, Poland\\
e-mail: mmroczko@mat.ug.edu.pl}

\begin{abstract}
We consider diagrams of links in $S^2$ obtained by projection from $S^3$ with the Hopf map
and the minimal crossing number for such diagrams. Knots admitting diagrams with at most one crossing are classified.
Some properties of these knots are exhibited. In particular, we establish which of these knots are algebraic and, for such knots, give an answer
to a problem posed by Fiedler in \cite{F}.
\end{abstract}
\maketitle

\let\thefootnote\relax 

\section{Introduction}
Let $p:S^3\to S^2$ be the Hopf map in the Hopf fibration of $S^3$.
For a link $L$ in $S^3$, we study the minimal number of crossings of $p(L)$ among generic projections in $S^2$.
Let us call this number the {\it Hopf crossing number} of $L$ and denote it by $h(L)$. The classical crossing number is denoted by $c(L)$. 

The Hopf crossing number was considered by Fiedler \cite{F} in the context of algebraic links coming from singularities of complex curves.
For an algebraic link $L$ realized as the intersection of a complex plane algebraic curve $X$ with a small sphere $S^3$ with a singularity of $X$
at its center, he considered the canonical Hopf fibration from $S^3$ to $\C P^1\approx S^2$, i.e. the fibration induced by the intersection of $S^3$ 
with the complex lines through the center of $S^3$. He defined $C_{alg}(L)$ as the minimal number of crossings under this Hopf fibration 
for all generic realizations of $L$ as $X\cap S^3$.
He computed $C_{alg}$ for torus knots $T(p,q)$, $1<p<q$, $q<2p$ and found lower bounds for $C_{alg}$ of other algebraic knots and links.

For an algebraic link $L$ one has $h(L)\le C_{alg}(L)$ as the link does not have to be realized as the intersection of a complex curve with a sphere.
In fact Fiedler poses a problem: is this inequality always an equality for algebraic links?

In this paper we classify knots with $h(K)\le 1$ (see Theorem~\ref{th_cl}) and show that, for such knots, 
the answer to Fiedler's problem is positive (see Theorem~\ref{th_fp}).
In order to achieve it, we use arrow diagrams introduced in \cite{MD}, which are equivalent to 
gleams introduced in \cite{T}.

In section~\ref{sec:arrow} we present the arrow diagrams of links in $S^3$ for which the minimal crossing number
is the Hopf crossing number.

In section~\ref{sec:examples} we consider several examples of arrow diagrams of knots.

In section~\ref{sec:jones} we compute the Jones polynomial of knots $K$ with $h(K)\le 1$.

In section~\ref{sec:classification} we classify knots $K$ with $h(K)\le 1$ and study some of their properties.
In particular, we find all knots $K$ with at most 10 crossings (in the classical sense) satisfying $h(K)\le 1$.

Finally, in section~\ref{sec:algknots} we establish which knots $K$, with $h(K)\le 1$, are algebraic and show that, for such knots, $h(K)=C_{alg}(K)$.

\section{Arrow diagrams of links}\label{sec:arrow}
\begin{figure}[h]
\centering
\includegraphics{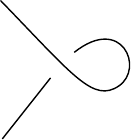} \raisebox{1.9 em}{$\; \overset{\Omega_1}{\longleftrightarrow} \;$} \includegraphics{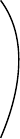}
\hspace{2em}
\includegraphics{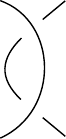} \raisebox{1.9 em}{$\; \overset{\Omega_2}{\longleftrightarrow} \;$} \includegraphics{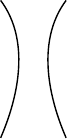}
\hspace{2em}
\includegraphics{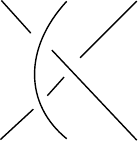} \raisebox{1.9 em}{$\; \overset{\Omega_3}{\longleftrightarrow} \;$} \includegraphics{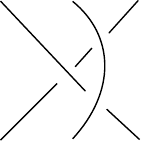}
\hspace{2em} \\[2em]
\includegraphics{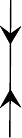} \raisebox{1.7 em}{$\; \overset{\Omega_4}{\longleftrightarrow} \;$} \includegraphics{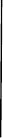}
\raisebox{1.7 em}{$\; \overset{\Omega_4}{\longleftrightarrow} \;$} \includegraphics{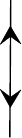} \hspace{2em}
\includegraphics{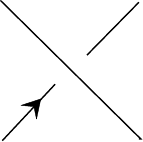} \raisebox{1.7 em}{$\; \overset{\Omega_5}{\longleftrightarrow} \;$} \includegraphics{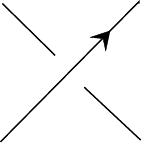}
\caption{Reidemeister moves}
\label{reid_moves}
\end{figure}

Let $T=S^1\times D$ be the solid torus. An {\it arrow diagram} of a link $L$ in $T$ is obtained by cuting
$T$ along $1\times D\times$, $1\in S^1$, and projecting $L$ from $[0,1]\times D$ thus obtained onto $D$. One keeps the information
about over- and undercrossings in this projection. Also, the points in $L\cap 1\times D$ are projected
onto {\it arrows} pointing to the part of $L$ that is close to $0$ in $[0,1]\times D$.
It was shown in~\cite{MD} that there are five Reidemeister moves connecting any two diagrams which represent
the same link in $T$: three classical ($\Om_1$, $\Om_2$ and $\Om_3$) 
and two extra ones ($\Om_4$ and $\Om_5$), presented in Figure~\ref{reid_moves}.

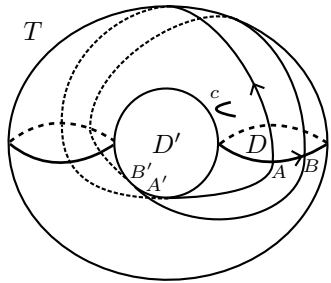
\begin{figure}\centering
\begin{overpic}{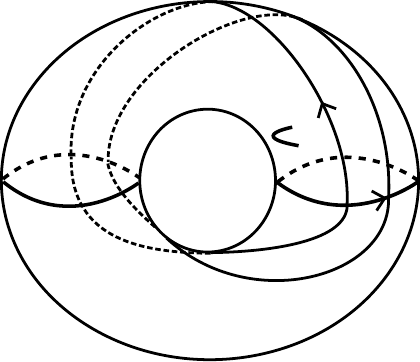}\put(92,34){\scriptsize $B$}\put(82,32){\scriptsize $A$}\put(74,40){$D$}\put(63,57){\scriptsize$c$}
\put(45,40){$D'$}\put(43,28){\scriptsize $A'$}\put(38,32){\scriptsize $B'$}\put(5,75){$T$}\end{overpic}
\caption{Hopf fibration}
\label{hopf_fibr}
\end{figure}

Given the solid torus $T$, $S^3$ is obtained from $T$ by gluing to it a second solid torus $T'$, see Figure~\ref{hopf_fibr}. 
In this figure, $T$ is shown. As a basis for $H_1(\partial T)$ we choose a standard longitude-meridian basis
of a torus embedded in a standard way in $S^3$.
In the decomposition $T=S^1\times D$, used for the arrow diagrams as above, we take as $S^1$ the curve passing through
$A$ and $A'$ (i.e. a fiber in the Hopf fibration of $S^3$). Such $S^1$, when oriented, represents $(1,1)$ in the chosen basis of $H_1(\partial T)$.
The solid torus $T'$ is glued to $T$ in such a way that its meridional disk $D'$, shown vertically in the middle
of the figure, is glued to $T$ along a $(1,0)$ curve in $\partial T$. Here, we choose the orientations of $S^1$ and $\partial D$
indicated on the figure by arrows.

\begin{figure}[h]
\includegraphics{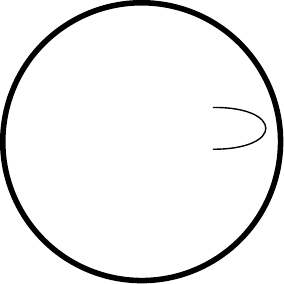} \raisebox{4.2 em}{$ \;
\underset{\longleftrightarrow}{\Omega_{\infty}} \; $ }
\begin{overpic}{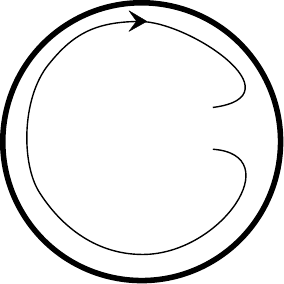}   \end{overpic}
\caption{The $\Om_{\infty}$ move}
\label{om_inf}
\end{figure}

Any link $L$ in $S^3$ can be pushed inside $T$. Its projection onto $D$ gives an arrow diagram of $L$,
as decribed above. If two links $L$ and $L'$, lying in $T$, are ambiently isotopic in $S^3$, one can pass from $L$ to $L'$ by a series
of ambient isotopies in $T$ and gliding of arcs along the disc $D'$. Such a gliding gives rise to one extra Reidemeister move, 
called the $\Om_{\infty}$ move, shown in Figure~\ref{om_inf}. Indeed, consider the small arc $c$ above $D$ in Figure~\ref{hopf_fibr}.
It projects onto an arc in $D$ as on the left of Figure~\ref{om_inf}. When $c$ is glided over $D'$ and slightly beyond, its projection in $D$
will be such as the one depicted on the right of Figure~\ref{om_inf}. Indeed, going around $\partial D$ in counterclockwise direction, one meets in succesion:
the unchanged endpoints of $c$, the arrow where the glided $c$ crosses $D$ next to $D\cap D'$ (the arrow is oriented in clockwise direction as we go opposite
to the orientation of the fibers $S^1$), a point next to $A'$ projecting to a point next to $A$ and, finally, a point next to $B'$ projecting to a point next to $B$.

Let $S^2$ be obtained from $D$ and $D'$ by identifying the points of $\partial D$ and $\partial D'$ belonging to the same fiber 
(such as $A$ and $A'$ or $B$ and $B'$ in Figure~\ref{hopf_fibr}). Then the Hopf map, $p:S^3\to S^2$, maps $T$ onto $D\subset S^2$ and $T'$ onto $D'\subset S^2$. 
As any link $L$ in $S^3$ can be pushed into $T$ without changing the number of crossings of $p(L)$, the Hopf crossing number $h(L)$ is the minimum of crossings
in all arrow diagrams of $L$.

For any link $L$ its classical diagrams form a subset of its arrows diagrams (consisting of diagrams
with no arrows). Thus, obviously, $h(L)\le c(L)$.

As in our later computations we use the Kauffman bracket, it is necessary to consider {\it framed} links. Such a link is equipped with a {\it framing} i.e.
a smooth family of unit vectors orthogonal to the link.
An arrow diagram represents a framed link where the framing is orthogonal to the fibers $S^1$ of the Hopf fibration.
All Reidemeister moves except $\Om_1$ preserve this framing. For $\Om_{\infty}$, notice that when we glide an arc, such as $c$ in Figure~\ref{hopf_fibr}, through
$D'$, the framing stays orthogonal to the fibers during the gliding.

Let $D$ be an arrow diagram. Denote by $\overline{w}(D)$ the writhe of the diagram obtained from $D$ by ignoring all arrows.
Denote by $w(D)$ the writhe of the framed link represented by $D$. If $D'$ has no arrows and is obtained from $D$ by Reidemeister moves
including $a$ moves $\Om_1$ increasing the writhe and $b$ moves $\Om_1$ decreasing the writhe, then $w(D)=w(D')+b-a$.

For an oriented diagram an arrow is {\it positive} if it points in the same direction as the orientation.
Otherwise it is {\it negative}. An arrow is {\it removable} if it can be eliminated with an $\Om_{\infty}$
move followed by an $\Om_4$ move. Thus, an arrow is removable if it is not separated from the boundary of
the diagram by any arc and is clockwise (so it will be opposite to the arrow created by $\Om_{\infty}$).

\begin{lemma}\label{lemma_writhe}
Suppose $D$ is an oriented arrow diagram of a knot with all arrows being removable. Suppose that there are $a>0$
positive and $b\ge 0$ negative arrows. Then $w(D)=\overline{w}(D)+a(a-1)+b(b-1)-2ab$.
\begin{proof}
Denote $D$ by $D_{a,b}$. Eliminate a positive arrow by $\Om_{\infty}$ and $\Om_4$, then push with $\Om_2$ and $\Om_5$
the other arrows through the new arc so that they are again removable. Denote the resulting diagram by
$D_{a-1,b}$. This is shown in Figure~\ref{pos_neg_arrows} and one sees that the remaining $a-1$ positive
arrows contribute $+2$ each to $\overline{w}(D_{a-1,b})$ whereas the $b$ negative arrows contribute $-2$ each.
Thus $\overline{w}(D_{a-1,b})=\overline{w}(D_{a,b})+2(a-1)-2b$.
We apply the same moves to $D_{a-1,b}$ getting
$\overline{w}(D_{a-2,b})=\overline{w}(D_{a-1,b})+2(a-2)-2b$.
Continuing until $a=0$ we get $\overline{w}(D_{0,b})=\overline{w}(D_{a,b})+2(a-1)-2b+2(a-2)-2b+...+2-2b+0-2b=\overline{w}(D_{a,b})+a(a-1)-2ab$.

Now reverse the orientation of $D_{0,b}$ to get $D_{b,0}$ ($\overline{w}$ is unchanged).
By the first part $\overline{w}(D_{0,0})=\overline{w}(D_{b,0})+b(b-1)$.
As $D_{0,0}$ has no arrows $\overline{w}(D_{0,0})=w(D_{0,0})$.

We have not used $\Om_1$ moves, so $w(D)=w(D_{a,b})=w(D_{0,0})$.
Thus $w(D)=\overline{w}(D)+a(a-1)-2ab+b(b-1)$.

\begin{figure}[h]
\centering
\includegraphics{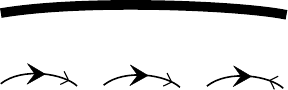} {$ \; \longrightarrow \; $ }
\includegraphics{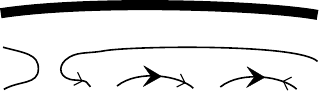} {$ \; \longrightarrow \; $ }
\includegraphics{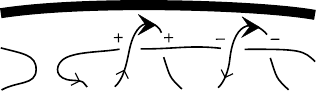} 
\caption{Removing a positive arrow}
\label{pos_neg_arrows}
\end{figure}

\end{proof}
\end{lemma}

We will later need to estimate, from an arrow diagram of a knot $K$, the classical crossing number $c(K)$. The following lemma will be useful:

\begin{lemma}\label{lemma_boundcross} Suppose that a knot $K$ has an arrow diagram $D$ with $k$ crossings, $a+b$ arrows of which $a$ are removable and $b$ are not removable
but are not separated from the boundary of the diagram.

Then $c(K)\le k+a(a-1)$ if $b=0$ and $c(K)\le k+b-1+(a+b)(a+b-1)$ if $b>0$.
\begin{proof}
If $b=0$, remove all $a$ arrows just as in the proof of the preceding lemma.
This creates $2(a-1)+2(a-2)+...+2=a(a-1)$ crossings and all arrows have been eliminated.

If $b>0$, create $b-1$ kinks with $\Om_1$ moves next to each arrow that is not removable except for one (so do nothing if $b=1$).
Then push such arrows into the kinks with $\Om_5$. This increases the number of crossings by $b-1$. Now $a+b-1$ arrows are removable. 
Remove these arrows as in the proof of the preceding lemma. This creates $2(a+b-2)+2(a+b-3)+...+2$ crossings. 
Push the last non removable arrow across the $a+b-1$ arcs with $\Om_2$ and $\Om_5$ creating $2(a+b-1)$ crossings.
At this stage there are $k+b-1+2(a+b-1)+2(a+b-2)+..+2=k+b-1+(a+b)(a+b-1)$ crossings.
Now create a kink, push the last arrow into it with $\Om_5$ and eliminate the arrow with $\Om_\infty$. As there are no more arrows, the arc created
by $\Om_\infty$ going next to the boundary of the diagram can be shrinked above the rest of the diagram and the crossing coming from the last kink
created can be eliminated with $\Om_1$. So we get the required number of crossings.
\end{proof}
\end{lemma}

\section{Examples: arrow diagrams of torus knots and some other knots}\label{sec:examples}
In this section we consider some arrow diagrams of torus knots and a couple of other knots that will be needed later.

Consider a torus knot $T(n,m)$, $n<m$, $gcd(n,m)=1$, lying on the boundary of the solid torus $T$ (see Figure~\ref{hopf_fibr}), 
going in a uniform way $n$ times in the longitudinal and $m$ times in the meridional direction. Thus, its homology class in $H_1(\partial T)$ is $(n,m)$,
in the standard longitude-meridian basis of $H_1(\partial T)$.

\begin{figure}\centering
\raisebox{1em}{\includegraphics{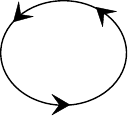}}\mbox{\qquad}\includegraphics{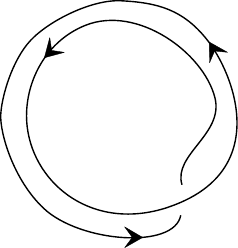}\mbox{\qquad}\includegraphics{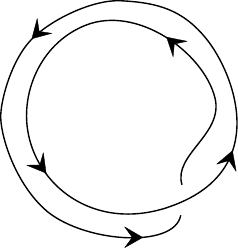}\mbox{\qquad}
\includegraphics{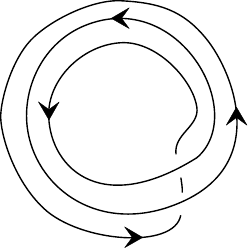}
\caption{Torus knots $T(3,4)$, $T(3,5)$, $T(5,7)$ and $T(4,7)$}
\label{ex_Txx}
\end{figure}

The Hopf map projects $T$ onto $D$ and $T(n,m)$ onto the circle $\partial D$.
Push the knot into $T$ and perturb it slightly, to get a generic projection with all crossings near each other (see Figure~\ref{ex_Txx}). 
We may choose the cut in $T$ (which creates the arrows) and the perturbation of the knot in such a way that there is an arrow just before all crossings and pointing
to the undercrossings. All other arrows are uniformly spaced on the resulting curve.
Now each arrow corresponds to going once along a Hopf fiber i.e. going around a curve $(1,1)$. Going around the disk corresponds to 
a meridional curve $(0,1)$. Thus, going $a$ times around the disk, with $b$ arrows, will correspond to a torus knot of type $a(0,1)+b(1,1)=(b+a,a)$
(the addition and multiplication in $H_1(\partial T)$). In Figure~\ref{ex_Txx} the first knot is of type $(1,1)+(1,1)+(1,1)+(0,1)=(3,4)$ i.e. the knot $T(3,4)$.
The second knot is of type $3(1,1)+2(0,1)=(3,5)$, the third is of type $5(1,1)+2(0,1)=(5,7)$ and the fourth is of type $4(1,1)+3(0,1)=(4,7)$.

It follows that an oval with $n$ counterclockwise arrows on it is the torus knot $T(n,n+1)$. Thus $h(T(n,n+1))=0$.
More generally we see that $T(n,n+k)$ has a diagram with $k-1$ crossings so $h(T(n,n+k))\le k-1$.
Indeed, consider a diagram $D$ going around the disk $k$ times with $k-1$ crossings, generalizing the diagrams of Figure~\ref{ex_Txx}
(for example $k=3$ for the last diagram in this figure).
Suppose that there are $n$ arrows on $D$ spaced uniformly, where the first arrow is just before the crossings. Then $D$ is a diagram of the
torus knot of type $k(0,1)+n(1,1)=(n,n+k)$.
Fiedler has shown in~\cite{F}, that $C_{alg}(T(n,n+k))=k-1$ if $k<n$. It does not imply that
$h(T(n,n+k))=k-1$ if $k<n$. We will get the last equality for $k\le 3$, see Proposition~\ref{prop_torus}.

If $k>n$, there are less arrows contributing $(1,1)$ than meridional curves contributing $(0,1)$ in the arrow diagrams such as in Figure~\ref{ex_Txx}.
Consider the same diagram $D$ with $k-1$ crossings as in the preceding paragraph: going around the disk $k$ times with $n$
arrows spaced uniformly, the first arrow just before the crossings. As $n<k$, the second arrow will not be on the most nested kink, so this
kink can be removed with $\Om_1$, yielding a diagram with $k-2$ crossings. Thus, $h(T(n,n+k))\le k-2$ if $k>n$. 
It is possible that more nested kinks may be removed. The most extreme case is $T(2,2+k)$ with a diagram containing only two arrows, 
so that all the most nested kinks up to the second arrow can be removed. One checks easily that $h(T(2,2+k))\le \frac{k-1}{2}$.

\begin{example}\label{ex:51}
Three arrow diagrams of the torus knot $T(2,5)$ are shown in Figure~\ref{ex_51}.

\begin{figure}\centering
\includegraphics{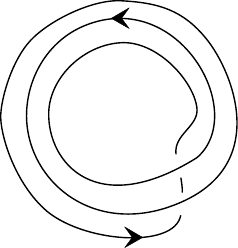}
\raisebox{1 em}{$ \; \overset{\Om_1}{\longrightarrow} \; $ }
\includegraphics{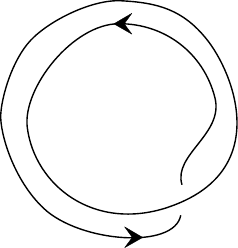}
\raisebox{1 em}{$ \; \overset{\Om_{\infty}}{\longrightarrow} \; $ }
\includegraphics{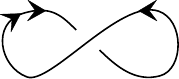}
\caption{Knot $T(2,5)$}
\label{ex_51}
\end{figure}

\end{example}

\begin{example}\label{ex:k21}
Consider the knot with a diagram shown in Figure~\ref{ex_942}. In this figure Reidemeister moves are used to show
that it is the knot $9_{42}$. In the last row of this figure arcs that are moved in the diagrams are thickened.
The third diagram from the end is transformed by a rotation by $\pi$ of a solid torus containing the knot indicated with
dashed lines.

\begin{figure}\centering
\includegraphics{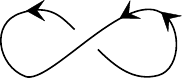}
\raisebox{1 em}{$ \; \overset{\Om_1, \Om_5}{\longrightarrow} \; $ }
\includegraphics{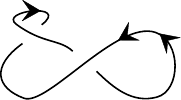}
\raisebox{1 em}{$ \; \overset{\Om_{\infty}, \Om_4}{\longrightarrow} \; $ }
\includegraphics{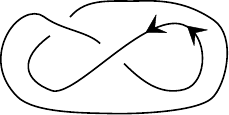}
\raisebox{1 em}{$ \; \overset{\Om_2, \Om_1, \Om_5}{\longrightarrow} \; $ } \\[2em]
\includegraphics{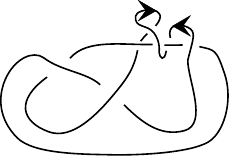}
\raisebox{1 em}{$ \; \overset{\Om_{\infty}, \Om_4}{\longrightarrow} \; $ }
\includegraphics{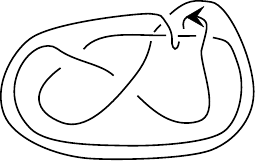}
\raisebox{1 em}{$ \; \overset{\Om_2,\Om_5,\Om_{\infty},\Om_4, \Om_1}{\longrightarrow} \; $ }
\includegraphics{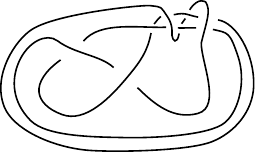} \\[2em]
\raisebox{1 em}{$\rightarrow \;\;$}
\includegraphics{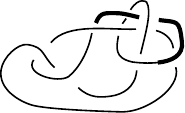}\raisebox{1 em}{$\;\; \rightarrow \;\;$}
\includegraphics{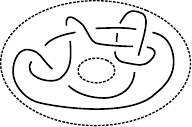}\raisebox{1 em}{$\;\; \rightarrow \;\;$}
\includegraphics{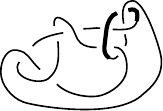}\raisebox{1 em}{$\;\; \rightarrow \;\;$}
\includegraphics{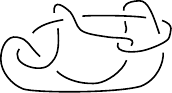}
\caption{Knot $9_{42}$}
\label{ex_942}
\end{figure}

\end{example}

\begin{example}\label{ex:k31}
Denote  the first arrow diagram in Figure~\ref{ex_942} by $K_{2,1}$. Denote a similiar diagram with 3 counterclockwise arrows on the right loop instead of 2
by $K_{3,1}$ (see Figure~\ref{K_relations} below). Later, we will see that $K_{2,1}$ and $K_{3,1}$ share the same Jones polynomial.
By eliminating all arrows the HOMFLY polynomial of these two knots can be computed with KnotPlot\cite{S}. One gets for $K_{2,1}$ and $K_{3,1}$ respectively:
\[P_{K_{2,1}}=-2l^{-2}-3-2l^2+m^2l^{-2}+4m^2+m^2l^2-m^4\]
\begin{align*}
P_{K_{3,1}}=5+10l^2+8l^4+2l^6-10m^2-25m^2l^2-14m^2l^4-m^2l^6\\
+6m^4+22m^4l^2+7m^4l^4-m^6-8m^6l^2-m^6l^4+m^8l^2
\end{align*}
Thus these two knots are distinct (including mirror images) and because of the assymetry of the $l$ terms $K_{3,1}$ is chiral.
One can check that the reductions of the two HOMFLY polynomials to the Jones polynomial give indeed the same polynomial.
\end{example}

\section{Jones polynomial of knots $K$ with $h(K)\le 1$}\label{sec:jones}
In this section we compute the Jones polynomial of knots with Hopf crossing number at most $1$.
This will allow us to classify such knots in the next section.

For $n\ge 0$, let $T_n$ be an oval with $n$ counterclockwise arrows on it. It is the right handed torus knot $T(n,n+1)$ (see the previous section).
The Jones polynomial of a torus knot $T(m,n)$ is given by \cite{J}:
\[V_{T(m,n)}=t^\frac{(m-1)(n-1)}{2}\frac{1-t^{m+1}-t^{n+1}+t^{m+n}}{1-t^2}\label{eq:VT1}\tag{VT1}\]
Thus, for right handed torus knots $T(n,n+1)$ the Jones polynomial is: 
\[V_{T_n}=t^\frac{n(n-1)}{2}\frac{1-t^{n+1}-t^{n+2}+t^{2n+1}}{1-t^2}\label{eq:VT2}\tag{VT2}\]
We extend to $n\in\Z$ the definition of $T_n$: if $n<0$, $T_n$ is an oval with $|n|$ clockwise arrows on it.
If $n<0$, $T_n$ can be transformed with $\Om_{\infty}$ and $\Om_4$ into $T_{-n-1}$.
In the same way, for $n\in\Z$, let $T'_n$ be an oval with $n$ clockwise arrows, if $n\ge 0$, and $|n|$ counterclockwise arrows, if $n<0$.
One can transform $T'_n$ into $T_{n-1}$ with $\Om_{\infty}$ and $\Om_4$, so these diagrams represent the same knot.
One checks easily that in the equation~(\ref{eq:VT2}) $V_{T_n}$ and $V_{T_{-n-1}}$ are equal for all $n$.
Thus, this equation is also valid for negative $n$.

An arrow diagram of a knot with no crossings is just an oval with some arrows on it. Opposite arrows can be removed with $\Om_4$.
Thus, from the preceding paragraph we get immediately (see also~\cite{F}):

\begin{lemma}\label{lemma_nocrossings}
A knot $K$ satisfies $h(K)=0$ if and only if $K$ is a torus knot $T(n,n+1)$.
\end{lemma}

We use the Kauffman bracket in order to compute the Jones polynomial (see \cite{K}).
Substituting $t=A^{-4}$, for a framed knot, $<K>=(-A)^{3w(K)}V(K)$,
where $<>$, $V$ and $w$ stand for the Kauffman bracket, Jones polynomial and writhe of the knot, respectively.
Until the end of this section we will use both $A$ and $t$, sometimes in the same formulas, with the understanding that $t=A^{-4}$.

Notice that the Kauffman bracket and Kauffman relations (smoothing a crossing and removing a trivial component)
are defined for the orthogonal projection of framed links from $\R^3$ to $\R^2$. However, the Kauffman is also defined for arrow diagrams, which, 
as classical diagrams, represent framed links in $S^3$. Also, Kauffman relations hold for arrow diagrams. Indeed, to apply a Kauffman relation to an arrow
diagram, we may first remove all arrows with appropriate Reidemeister moves, apply the relation to a classical diagram, then recreate the removed
arrows with Reidemeister moves. 

\begin{lemma}
The writhe of $T_n$ is $n(n+1)$, for $n\in\Z$.
\begin{proof}
If $n=0$, then it is trivial. If $n>0$, then apply $\Om_{\infty}$ to $T_n$ to get $T'_{n+1}$. This does not change the writhe. 
Orient the diagram so that all the arrows of $T'_{n+1}$ are positive and removable as in the definition before lemma~\ref{lemma_writhe}. 
From this lemma, $w(T'_{n+1})=\overline{w}(T'_{n+1})+(n+1)n=n(n+1)$.

If $n<0$, orient the diagram so that there are $|n|$ postive removable arrows. From the same lemma,
$w(T_n)=w(T'_{|n|})=|n|(|n|-1)=(-n)(-n-1)=n(n+1)$.
\end{proof}
\end{lemma}

Using equation~(\ref{eq:VT2}) and the last lemma, one gets: 
\[<T_n>=A^{3n(n+1)}t^\frac{n(n-1)}{2}\frac{1-t^{n+2}-t^{n+1}(1-t^n)}{1-t^2}\]

Let us write $<T_n>=U_n+V_n$ where $U_n$ contains $1-t^{n+2}$ and $V_n$
contains $-t^{n+1}(1-t^n)$ of the numerator of the bracket.

Let us also write $<T_n>=U'_n+V'_n$ where $U'_n$ contains $1-t^{n+1}$ and $V'_n$
contains $-t^{n+2}(1-t^{n-1})$ of the numerator of the bracket.

\begin{lemma}\label{lemma_UV}
Let $n\in\Z$. Then $V_n+A^2U_{n-2}=0$ and $V'_n+A^{-2}U'_{n-2}=0$.
\begin{proof}
\[V_n=-t^{n+1} A^{3n(n+1)} t^\frac{n(n-1)}{2}\frac{1-t^n}{1-t^2}=-A^v\frac{1-t^n}{1-t^2}\]
\[A^2U_{n-2}=A^2 A^{3(n-2)(n-1)}t^\frac{(n-2)(n-3)}{2}\frac{1-t^n}{1-t^2}=A^u\frac{1-t^n}{1-t^2}\]
One has to check that $v=u$ so these two terms cancel each other.
\[v=-4(n+1)+3n(n+1)-2n(n-1)=n^2+n-4\]
\[u=2+3(n-2)(n-1)-2(n-2)(n-3)=n^2+n-4\]
The proof for $V'_n$ and $U'_{n-2}$ is similar. In fact, using the same $v$ and $u$ as above one has:
\[V'_n=-tA^v\frac{1-t^{n-1}}{1-t^2}\]
\[A^{-2}U'_{n-2}=tA^2U'_{n-2}=tA^u\frac{1-t^{n-1}}{1-t^2}\]
As $v=u$, the equality follows.
\end{proof}
\end{lemma}

\begin{lemma}\label{lemma_S}
For $n,m\in\Z$, $m>0$ let:
\begin{eqnarray*}
S_{n,m} &=& <T_n>+A^2<T_{n-2}>+...+A^{2m}<T_{n-2m}>\\
S'_{n,m}&=&<T_n>+A^{-2}<T_{n-2}>+...+A^{-2m}<T_{n-2m}>
\end{eqnarray*}
Then:
\begin{eqnarray*}
S_{n,m} &=&A^{3n(n+1)}t^\frac{n(n-1)}{2}\left[\frac{1-t^{n+2}-t^{(m+1)(n-m+1)}(1-t^{n-2m})}{1-t^2}\right]\\
S'_{n,m}&=&A^{3n(n+1)}t^\frac{n(n-1)}{2}\left[\frac{1-t^{n+1}-t^{(m+1)(n-m+2)}(1-t^{n-2m-1})}{1-t^2}\right]
\end{eqnarray*}
\begin{proof}
From Lemma~\ref{lemma_UV}, $S_{n,m}=U_n+A^{2m}V_{n-2m}$, as the other terms cancel each other. Thus:
\begin{eqnarray*}
S_{n,m}&=&A^{3n(n+1)}t^\frac{n(n-1)}{2}\frac{1-t^{n+2}}{1-t^2}\\
&&+A^{2m}A^{3(n-2m)(n-2m+1)}t^\frac{(n-2m)(n-2m-1)}{2}\left(-t^{n-2m+1}\frac{1-t^{n-2m}}{1-t^2}\right)
\end{eqnarray*}
For the second term:
\begin{eqnarray*}
&&A^{2m} A^{3(n-2m)(n-2m+1)}=A^{3n(n+1)}A^{2m} A^{-6m(n-2m+1)-6mn}\\
&=&A^{3n(n+1)}A^{-12mn+12m^2-4m}=A^{3n(n+1)}t^{3mn-3m^2+m}
\end{eqnarray*}
and
\begin{eqnarray*}
&&t^\frac{(n-2m)(n-2m-1)}{2}t^{n-2m+1}=t^\frac{n(n-1)}{2}t^{-m(n-2m-1)-mn+n-2m+1}\\
&=&t^\frac{n(n-1)}{2}t^{-2mn+2m^2-m+n+1}
\end{eqnarray*}
Thus, factoring out $A^{3n(n+1)}t^\frac{n(n-1)}{2}$ in the second term one gets:
\[-t^{-2mn+2m^2-m+n+1+3mn-3m^2+m}=-t^{mn-m^2+n+1}=-t^{(m+1)(n-m+1)}\]
For the second part, again from Lemma~\ref{lemma_UV}, $S'_{n,m}=U'_n+A^{-2m}V'_{n-2m}$. Thus:
\begin{eqnarray*}
S'_{n,m}&=&A^{3n(n+1)}t^\frac{n(n-1)}{2}\frac{1-t^{n+1}}{1-t^2}\\
&&+A^{-2m}A^{3(n-2m)(n-2m+1)}t^\frac{(n-2m)(n-2m-1)}{2}\left(-t^{n-2m+2}\frac{1-t^{n-2m-1}}{1-t^2}\right)
\end{eqnarray*}
The difference between the coefficient of the second term in $S'_{n,m}$ compared to $S_{n,m}$ is the multplication by $A^{-2m}$ instead of $A^{2m}$
(which corresponds to multiplying by $A^{-4m}=t^m$) and $-t^{n-2m+2}$ instead of $-t^{n-2m+1}$. Overall, this corresponds
to multiplying by $t^{m+1}$.
Thus, factoring out $A^{3n(n+1)}t^\frac{n(n-1)}{2}$ in the second term one gets:
\[-t^{(m+1)(n-m+1)}t^{m+1}=-t^{(m+1)(n-m+2)}\]
\end{proof}
\end{lemma}

From the definition of the Kauffman bracket (see \cite{K}) one gets easily the two following relations:

\begin{align*}
\centering
\raisebox{-7pt}{\includegraphics{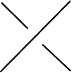}}=A^2 \raisebox{-7pt}{\includegraphics{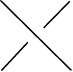}}
+(A^{-1}-A^{3}) \raisebox{-7pt}{\includegraphics{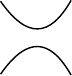}}
\end{align*}

\begin{align*}
\centering
\raisebox{-7pt}{\includegraphics{Lminus}}=A^{-2} \raisebox{-7pt}{\includegraphics{Lplus}}
+(A-A^{-3}) \raisebox{-7pt}{\includegraphics{Lzero}}
\end{align*}

For $a\ge 0$, $b\ge 0$, let $K_{a,b}$, $K'_{a,b}$ and $K''_{a,b}$ be the framed knots shown on the left of Figure~\ref{K_relations}, where $a$ next
to an arrow stands for $a$ arrows (and the same for $b$).
In this figure the Kauffman bracket of these knots is calculated using the previous two relations (for simplicity, the bracket is omitted in the figure). 
For example, applying $\Om_5$ and $\Om_4$ moves, in the first case one gets:

$<K_{a,b}>=A^2<K_{a-1,b-1}>+(A^{-1}-A^3)<T_{a+b}>$ 

Iterating the first relation until $b=0$ and removing the positive kink with $-A^{3}$ in the Kauffman bracket one gets:

\begin{align*}
\tag{1}\label{eq:Kab}<K_{a,b}>=(A^{-1}-A^3)[<T_{a+b}>+A^2 <T_{a+b-2}>+\ldots\\
\ldots+A^{2b-2}<T_{a-b+2}>]-A^{2b+3}<T_{a-b}>
\end{align*}

Symmetrically, exchanging $A$ and $A^{-1}$ and $T_k$ with $T'_k=T_{k-1}$, one gets for $K'$:

\begin{align*}
\label{eq:Kpab}\tag{2}<K'_{a,b}>=(A-A^{-3})[<T_{a+b-1}>+A^{-2} <T_{a+b-3}>+\ldots\\
\ldots+A^{-2b+2}<T_{a-b+1}>]-A^{-2b-3}<T_{a-b-1}>
\end{align*}

For $K''$ one gets:

\begin{align*}
\label{eq:Kppab}\tag{3}<K''_{a,b}>=(A^{-1}-A^3)[<T_{a-b-1}>+A^2 <T_{a-b+1}>+\ldots\\
\ldots+A^{2b-2}<T_{a+b-3}>]-A^{2b+3}<T_{a+b-1}>
\end{align*}

\begin{figure}\centering
\raisebox{10pt}{$K_{a,b}:=\;\;$}
\oversimage{revK_1}{0,47}{\scriptsize $b$}{80,47}{\scriptsize $a$}
\raisebox{10pt}{$\;\;=\;\;A^2\;\;$}
\oversimage{revK_1}{-10,47}{\scriptsize $b-1$}{75,47}{\scriptsize $a-1$}
\raisebox{10pt}{$\;+(A^{-1}-A^3)\;\;$}
\overimage{revert_3b}{-10,47}{\scriptsize $a+b$} \\[2em]
\raisebox{10pt}{$K'_{a,b}:=\;\;$}
\oversimage{rev1}{0,47}{\scriptsize $a$}{80,47}{\scriptsize $b$}
\raisebox{10pt}{$\;\;=\;\;A^{-2}\;\;$}
\oversimage{rev1}{-10,47}{\scriptsize $a-1$}{75,47}{\scriptsize $b-1$}
\raisebox{10pt}{$\;\;+(A-A^{-3})\;\;$}
\overimage{revert_3}{-10,47}{\scriptsize $a+b$} \\[2em]
\raisebox{10pt}{$K''_{a,b}:=\;\;$}
\oversimage{revert_1}{0,47}{\scriptsize $a$}{80,47}{\scriptsize $b$}
\raisebox{10pt}{$\;\;=\;\;A^2\;\;$}
\oversimage{revert_1}{-10,47}{\scriptsize $a+1$}{75,47}{\scriptsize $b-1$}
\raisebox{10pt}{$\;\;+(A^{-1}-A^3)\;\;$}
\overimage{revert_3}{-10,47}{\scriptsize $a-b$}
\caption{Relations for K, K' and K''}
\label{K_relations}
\end{figure}

\begin{proposition}\label{prop_bracket}
Let $n=a+b$, then: 
\begin{eqnarray*}
<K_{a,b}>&=&A^{3n(n+1)}t^\frac{n(n-1)}{2}A^{-1}\left[\frac{1-t^{n+2}-t^{b(a+2)}(1-t^{a-b+2})}{1-t^2}\right.\\
& &\left.-t^{-1}\frac{1-t^{n+2}-t^{(b+1)(a+1)}(1-t^{a-b})}{1-t^2}\right]\\
<K'_{a,b}>&=&A^{3n(n-1)}t^\frac{(n-1)(n-2)}{2}A\left[\frac{1-t^n-t^{b(a+2)}(1-t^{a-b})}{1-t^2}\right.\\
&&\left.-t\frac{1-t^n-t^{(b+1)(a+1)}(1-t^{a-b-2})}{1-t^2}\right]\\
<K''_{a,b}>&=&A^{3n(n-1)}A^{2b+3}t^\frac{(n-1)(n-2)}{2}\left[t^{n+2}\frac{1-t^{n-2}-t^{ba}(1-t^{a-b-2})}{1-t^2}\right.\\
&&\left.-\frac{1-t^n-t^{(b+1)(a+1)}(1-t^{a-b-2})}{1-t^2}\right]
\end{eqnarray*}

\begin{proof}
Using $S_{n,m}$ and $S'_{n,m}$ as in Lemma~\ref{lemma_S}, equation~(\ref{eq:Kab}) can be rewritten as:
\[K_{a,b}=A^{-1}S_{n,b-1}-A^3S_{n,b}=A^{-1}(S_{n,b-1}-t^{-1}S_{n,b})\]
Similarily, equation~(\ref{eq:Kpab}) can be rewritten as:
\[K'_{a,b}=AS'_{n-1,b-1}-A^{-3}S'_{n-1,b}=A(S'_{n-1,b-1}-tS'_{n-1,b})\]
Equation~(\ref{eq:Kppab}) can be rewritten as:
\[K''_{a,b}=A^{-1}A^{2b-2}S'_{n-3,b-1}-A^3A^{2b}S'_{n-1,b}=A^{2b-3}S'_{n-3,b-1}-A^{2b+3}S'_{n-1,b}\]
The result follows from the formulas for $S$ and $S'$ from Lemma~\ref{lemma_S}. For the first term of $<K''_{a,b}>$ one uses
the identity that is easily checked:
\[A^{2b-3}A^{3(n-3)(n-2)}t^\frac{(n-3)(n-4)}{2}=A^{3n(n-1)}A^{2b+3}t^\frac{(n-1)(n-2)}{2}t^{n+2}\]
\end{proof}
\end{proposition}

To get the Jones polynomial from the Kauffman bracket it suffices to compute the writhes of the knots.

\begin{lemma}\label{lemma_wrK}
$w(K_{a,b})=1+a(a+1)+b(b+1)-2ab$\\
$w(K'_{a,b})=-1+a(a-1)+b(b-1)-2ab$\\
$w(K''_{a,b})=1+a(a-1)+b(b+1)+2ab$
\begin{proof}
In order to use lemma~\ref{lemma_writhe} we want all arrows to be removable. This is the case for $K'_{a,b}$.
The single crossing of $K'_{a,b}$ is negative. Also we can orient $K'_{a,b}$ in such a way that the $a$ arrows on the
left loop are positive and the $b$ arrows on the right loop are negative.
From lemma~\ref{lemma_writhe} we get:

$w(K'_{a,b})=\overline{w}(K'_{a,b})+a(a-1)+b(b-1)-2ab=-1+a(a-1)+b(b-1)-2ab$

In $K_{a,b}$ none of the arrows are removable. Let $D_{a,b}$ be the diagram obtained from $K_{a,b}$ by putting with $\Om_1$ a negative kink
next to each arrow then pushing each arrow with $\Om_5$ into each such kink which become positive.
The $\Om_1$ moves change the writhe so that $w(K_{a,b})=w(D_{a,b})+a+b$. Because of the positive kinks $\overline{w}(D_{a,b})=\overline{w}(K_{a,b})+a+b=1+a+b$,
as the single crossing of $K_{a,b}$ is positive.
Now in $D_{a,b}$ all arrows are removable and it can be oriented so that there are $a$ postive and $b$ negative arrows.
Applying again lemma~\ref{lemma_writhe} we get:

$w(K_{a,b})=w(D_{a,b})+a+b=\overline{w}(D_{a,b})+a(a-1)+b(b-1)-2ab+a+b=a(a-1)+b(b-1)-2ab+a+b+a+b+1=1+a(a+1)+b(b+1)-2ab$

In $K''_{a,b}$ the $a$ arrows on the left kink are removable and the $b$ arrows on the right kink are not. As in the previous case, we make these $b$
arrows removable by $\Om_1$ and $\Om_5$ moves getting a diagram $D''_{a,b}$ with $w(K''_{a,b})=w(D''_{a,b})+b$ and 
$\overline{w}(D''_{a,b})=\overline{w}(K''_{a,b})+b=1+b$, as the single crossing of $K''_{a,b}$ is positive.
$D''_{a,b}$ can be oriented so that it has $a+b$ positive removable arrows and from lemma~\ref{lemma_writhe} we get:

$w(K''_{a,b})=w(D''_{a,b})+b=\overline{w}(D''_{a,b})+(a+b)(a+b-1)+b=(a+b)(a+b-1)+2b+1=
1+a(a-1)+ab+b(a+b-1)+2b+1=1+a(a-1)+2ab+b(b+1)$
\end{proof}
\end{lemma}

\begin{theorem}\label{th_jones}
Let $n=a+b$, then:
\begin{eqnarray*}
V_{K_{a,b}}&=&t^\frac{a^2+b^2-4ab-a-b}{2}\left[\frac{1-t^{n+2}-t^{(b+1)(a+1)}(1-t^{a-b})}{1-t^2}\right.\\
& &\left.-t\frac{1-t^{n+2}-t^{b(a+2)}(1-t^{a-b+2})}{1-t^2}\right]\\
V_{K'_{a,b}}&=&t^\frac{a^2+b^2-4ab-3a-3b}{2}\left[t\frac{1-t^n-t^{(b+1)(a+1)}(1-t^{a-b-2})}{1-t^2}\right.\\
&&\left.-\frac{1-t^n-t^{b(a+2)}(1-t^{a-b})}{1-t^2}\right]\\
V_{K''_{a,b}}&=&t^b t^\frac{(n-1)(n-2)}{2}\left[-t^{n+2}\frac{1-t^{n-2}-t^{ba}(1-t^{a-b-2})}{1-t^2}\right.\\
&&\left.+\frac{1-t^n-t^{(b+1)(a+1)}(1-t^{a-b-2})}{1-t^2}\right]
\end{eqnarray*}
\begin{proof}
As $V_K=(-A)^{-3w(K)}<K>$ let us compute the contributions of the writhe from Lemma~\ref{lemma_wrK} 
times the coefficients in front of the parentheses for
the Kauffman bracket of $K_{a,b}$, $K'_{a,b}$ and $K''_{a,b}$ from Proposition~\ref{prop_bracket}.

For $V_{K_{a,b}}$ it is:
\[(-A)^{-3(1+a(a+1)+b(b+1)-2ab)}A^{3n(n+1)}t^\frac{n(n-1)}{2}A^{-1}=\]
\[-A^{-3-3a^2-3a-3b^2-3b+6ab+3a^2+3b^2+6ab+3a+3b-2a^2-2b^2-4ab+2a+2b-1}=\]
\[-A^{-4-2a^2-2b^2+8ab+2a+2b}=-tt^\frac{a^2+b^2-4ab-a-b}{2}\]
Multpliying the terms in the parenthesis of $<K_{a,b}>$ by $-t$ we get the required Jones polynomial.

For $V_{K'_{a,b}}$ one has:
\[(-A)^{-3(-1+a(a-1)+b(b-1)-2ab)}A^{3n(n-1)}t^\frac{(n-1)(n-2)}{2}A=\]
\[-A^{3-3a^2+3a-3b^2+3b+6ab+3a^2+3b^2+6ab-3a-3b-2a^2-2b^2-4ab+6a+6b-4+1}=\]
\[-A^{-2a^2-2b^2+8ab+6a+6b}=-t^\frac{a^2+b^2-4ab-3a-3b}{2}\]

For $V_{K''_{a,b}}$ one has:
\[(-A)^{-3(1+a(a-1)+b(b+1)+2ab)}A^{3n(n-1)}A^{2b+3}t^\frac{(n-1)(n-2)}{2}=\]
\[-A^{-3-3a^2+3a-3b^2-3b-6ab+3a^2+3b^2+6ab-3a-3b+2b+3-2a^2-2b^2-4ab+6a+6b-4}=\]
\[-A^{-4-2a^2-2b^2-4ab+6a+2b}=-A^{-4b-2(a+b-1)(a+b-2)}=-t^bt^\frac{(n-1)(n-2)}{2}\]

\end{proof}
\end{theorem}

\section{Classification of knots $K$ with $h(K)\le 1$}\label{sec:classification}

\begin{lemma}\label{lemma_reducetoK}
Let $K$ be a knot with $h(K)=1$. Then $K$ is one of $K_{a,b}$, $K'_{a,b}$ with $a\ge b\ge 1$ or
$K''_{a,b}$ with $a>1$ and $b\ge 1$.
\begin{proof}
By assumption $K$ has an arrow diagram with one crossing i.e. with two loops. If one loop is nested in the other, perform an $\Om_{\infty}$
move to get two unnested loops (so the diagram now looks like $\infty$ with a positive or negative crossing). Let us call the loops
of the diagram of $K$ a {\it left} loop and a {\it right} loop.
On each loop, we eliminate opposite arrows with $\Om_4$. Now, if there are no arrows on one of the loops it can be eliminated with $\Om_1$ and
$h(K)=0$, a contradiction. Thus, there are some arrows on each loop.

Suppose first that the arrows on each loop are counterclockwise (as in $K_{a,b}$). If the unique crossing is different from $K_{a,b}$,
then we push one arrow from the left loop to the right with $\Om_5$ changing the crossing. Then we eliminate opposite arrows with $\Om_4$.
There should be arrows on each loop otherwise, as above, $h(K)=0$. If there are more arrows on the left loop than on the right loop, then
rotate the diagram of $K$ by $\pi$. The crossing is unchanged, the arrows are counterclockwise and we get $K_{a,b}$ with $a\ge b\ge 1$.

Suppose now that the arrows on each loop are clockwise (as in $K'_{a,b}$). By the same method as above we can change the crossing if needed
and put more arrows on the left by rotating the diagram. Thus we get $K'_{a,b}$ with $a\ge b\ge 1$.

Finally suppose that the arrows are counterclockwise on one loop and clockwise on the other (as in $K''_{a,b}$). First, modify as above the crossing
if needed. If the counterclockwise arrows are on the left loop, then rotate the diagram by $\pi$. In this way we get $K''_{a,b}$. If there is only
one arrow on the left loop it can be pushed to the right loop and $h(K)=0$. Thus we may assume that $a>1$. Notice that in this case it
is not necessary that $a\ge b$.
\end{proof}
\end{lemma}

We will use the Jones polynomials from Theorem~\ref{th_jones} to distinguish the knots $K_{a,b}$, $K'_{a,b}$ and $K''_{a,b}$.
First, one can consider their Jones polynomials multiplied by $1-t^2$ and by a power of $t$ so that 
the lowest coefficient has degree zero. Denote by $J_K$ such polynomial obtained from the Jones polynomial $V_K$ of a
knot $K$. For $K_{a,b}$, $a\ge b\ge 1$ one has:
\[ J_{K_{a,b}}=1-t-t^{a+b+2}+t^{a+b+3}+t^{ab+2b+1}-t^{ab+a+b+1}-t^{ab+a+b+3}+t^{ab+2a+1} \]
In this formula, we want to have increasing powers of $t$ and check for possible cancellations. One has to consider some
special cases.

First, suppose that $b=1$. Then:
\[ J_{K_{a,1}}=1-t+t^{a+4}-t^{2a+2}-t^{2a+4}+t^{3a+1} \label{eq:JKb1}\tag{JK1} \]
In the last forumla the powers are increasing except for
$a=1$ giving $J_{K_{1,1}}=1-t+t^5-t^6$ and $a=2$ or $a=3$ 
giving $J_{K_{2,1}}=J_{K_{3,1}}=1-t+t^7-t^8$.

Suppose now that $b>1$.
If $a=b$, $a=b+1$ or $a=b+2$, one has, respectively,
\[ J_{K_{b,b}}=1-t-t^{2b+2}+t^{2b+3}+t^{b^2+2b+1}-t^{b^2+2b+3}  \label{eq:JKaeqb}\tag{JK2} \]
\[ J_{K_{b+1,b}}=1-t-t^{2b+3}+t^{2b+4}+t^{b^2+3b+1}-t^{b^2+3b+2}+t^{b^2+3b+3}-t^{b^2+3b+4} \label{eq:JKaeqb1}\tag{JK3} \]
\[ J_{K_{b+2,b}}=1-t-t^{2b+4}+t^{2b+5}+t^{b^2+4b+1}-t^{b^2+4b+3} \label{eq:JKaeqb2}\tag{JK4} \]
Finally, if $a>b+2$,
\[ J_{K_{a,b}}=1-t-t^{a+b+2}+t^{a+b+3}+t^{ab+2b+1}-t^{ab+a+b+1}-t^{ab+a+b+3}+t^{ab+2a+1} \label{eq:JKagb2}\tag{JK5} \]

Similarily, for $K'_{a,b}$, $a\ge b\ge 1$ one has:
\[ J_{K'_{a,b}}=-1+t+t^{a+b}-t^{a+b+1}+t^{ab+2b}-t^{ab+a+b}-t^{ab+a+b+2}+t^{ab+2a} \]

Again, we consider special cases.
Suppose first that $b=1$. Then one has:
\[ J_{K'_{a,1}}=-1+t+t^{a+1}-t^{2a+1}-t^{2a+3}+t^{3a} \label{eq:JKpb1}\tag{JK'1} \]
In the last formula the powers are increasing except for
$a=1$ giving $J_{K'_{1,1}}=-1+t+t^2-t^5$, $a=2$ 
giving $J_{K'_{2,1}}=-1+t+t^3-t^5+t^6-t^7$
and $a=3$ giving $J_{K'_{3,1}}=-1+t+t^4-t^7$.

Suppose now that $b>1$.
If $a=b$, $a=b+1$ or $a=b+2$, one has, respectively,
\[ J_{K'_{b,b}}=-1+t+t^{2b}-t^{2b+1}+t^{b^2+2b}-t^{b^2+2b+2} \label{eq:JKpaeqb}\tag{JK'2} \]
\[ J_{K'_{b+1,b}}=-1+t+t^{2b+1}-t^{2b+2}+t^{b^2+3b}-t^{b^2+3b+1}+t^{b^2+3b+2}-t^{b^2+3b+3} \label{eq:JKpaeqb1}\tag{JK'3} \]
\[J_{K'_{b+2,b}}=-1+t+t^{2b+2}-t^{2b+3}+t^{b^2+4b}-t^{b^2+4b+2}  \label{eq:JKpaeqb2}\tag{JK'4} \]
Finally, if $a>b+2$,
\[J_{K'_{a,b}}=-1+t+t^{a+b}-t^{a+b+1}+t^{ab+2b}-t^{ab+a+b}-t^{ab+a+b+2}+t^{ab+2a}  \label{eq:JKpagb}\tag{JK'5} \]

Consider the third family $K''_{a,b}$ now. We suppose that $a>1$, $b\ge 1$ but it is not
necessary that $a\ge b$. One has:
\[ J_{K''_{a,b}}=1-t^{a+b}-t^{a+b+2}+t^{2a+2b}-t^{ab+a+b+1}+t^{ab+a+b+2}+t^{ab+2a-1}-t^{ab+2a} \]
Suppose first that $b=1$. One has:
\[ J_{K''_{a,1}}=1-t^{a+1}-t^{a+3}+t^{2a+3}+t^{3a-1}-t^{3a} \label{eq:JKppb1}\tag{JK''1} \]
In the last formula the powers are increasing except for
$a=2$ giving $J_{K''_{2,1}}=1-t^3-t^6+t^7$,
$a=3$ giving $J_{K''_{3,1}}=1-t^4-t^6+t^8$
and $a=4$ giving $J_{K''_{4,1}}=1-t^5-t^7+2t^{11}-t^{12}$.

Suppose now that $b>1$. As a special case suppose that $a=2$, then,
\[ J_{K''_{2,b}}=1-t^{b+2}-t^{b+4}+t^{2b+3}-t^{3b+3}+t^{3b+4}  \label{eq:JKppa2}\tag{JK''2} \]
Now suppose that $a>2$ (and $b>1$).
There may be some cancellations in the four last terms of the general formula. One has:
$(ab+2a-1)-(ab+a+b+1)=a-b-2$. So there are some terms with the same degree if $a-b$ is $1$, $2$ or $3$.

If $a\le b$,
\[ J_{K''_{a,b}}=1-t^{a+b}-t^{a+b+2}+t^{2a+2b}+t^{ab+2a-1}-t^{ab+2a}-t^{ab+a+b+1}+t^{ab+a+b+2} \label{eq:JKppalb}\tag{JK''3} \]
If $a=b+1$, $a=b+2$ or $a=b+3$, one has, respectively,
\[ J_{K''_{b+1,b}}=1-t^{2b+1}-t^{2b+3}+t^{4b+2}+t^{b^2+3b+1}-2t^{b^2+3b+2}+t^{b^2+3b+3} \label{eq:JKppaeqb1}\tag{JK''4} \] 
\[ J_{K''_{b+2,b}}=1-t^{2b+2}-t^{2b+4}+t^{4b+4} \label{eq:JKppaeqb2}\tag{JK''5} \]
\[ J_{K''_{b+3,b}}=1-t^{2b+3}-t^{2b+5}+t^{4b+6}-t^{b^2+5b+4}+2t^{b^2+5b+5}-t^{b^2+5b+6}\label{eq:JKppaeqb3}\tag{JK''6} \]
Finally, if $a>b+3$,
\[J_{K''_{a,b}}=1-t^{a+b}-t^{a+b+2}+t^{2a+2b}-t^{ab+a+b+1}+t^{ab+a+b+2}+t^{ab+2a-1}-t^{ab+2a} \label{eq:JKppagb}\tag{JK''7} \]

\begin{theorem}\label{th_cl}
The knots $K_{a,b}$, $K'_{a,b}$, $a\ge b\ge 1$ and $K''_{a,b}$, $a>1$, $b\ge 1$ are distinct.
Any knot $K$ with $h(K)=1$ is such a knot.

\begin{proof}
If two such knots are the same then their polynomials $J$ must be the same, so it is sufficient to check
that all these polynomials are different (except for $K_{2,1}$ and $K_{3,1}$ which are distinguished by the HOMFLY polynomial,
see example~\ref{ex:k31} in section~\ref{sec:examples}).

First we notice that $J_{K'_{a,b}}$ starts with $-1$ whereas $J_{K_{a,b}}$ and $J_{K''_{a,b}}$ start with $1$.
Also, the second term in $K_{a,b}$ is $-t$ whereas in $K''_{a,b}$ it is $-t^{a+b}$ and $a+b\ge 3$.
Thus these three families of knots are distinguished by $J$.

$\bullet$ Consider the family $K_{a,b}$. In the five formulas for $J$ above there are $6$ or $8$ terms except for the
special cases $K_{1,1}$, $K_{2,1}$ and $K_{3,1}$ for which there are $4$ terms. Thus, these special cases are
distinct from the rest (and $K_{1,1}$ is distinct from $K_{1,2}$ and $K_{1,3}$).

Consider the formulas with six terms i.e. equations (\ref{eq:JKb1}), (\ref{eq:JKaeqb}) and (\ref{eq:JKaeqb2}).
In equation (\ref{eq:JKb1}) the last term has coefficient $1$, which distinguishes it from the other two equations
with last coefficient $-1$.
For these two equations, suppose that $J_{K_{b,b}}=J_{K_{c+2,c}}$ for some integers $b>1$ and $c>1$.
By comparing the different terms in these polynomials we must have:

$2b+2=2c+4$ and $b^2+2b+1=c^2+4c+1$. From the first equation $b=c+1$ so $b^2+2b+1=c^2+4c+4\neq c^2+4c+1$.

Also, the knots $K_{a,1}$ are distinguished by $a$, $K_{b,b}$ are distinguished by $b$ and $K_{b+2,b}$ are
also distinguished by $b$.

Now consider the formulas with eight terms i.e. equations (\ref{eq:JKaeqb1}) and (\ref{eq:JKagb2}).
These two subfamilies are distinguished by the coefficient of the last term ($-1$ and $1$).
The knots $K_{b+1,b}$ are distinguished by $b$.
For the general $J_{K_{a,b}}$ of equation (\ref{eq:JKagb2}), notice that the difference between the powers of the second
and third term determines $a+b$, whereas the difference between the powers of the last two terms determines $a-b$.
Together, this determines the couple $(a,b)$ so all such knots are distinct.

$\bullet$ Consider the family $K'_{a,b}$. In the five formulas for $J$ above there are $6$ or $8$ terms except for the
specials cases $K'_{1,1}$, $K'_{3,1}$ for which there are $4$ terms. Thus, these special cases are
distinct from the rest (and $K'_{1,1}$ is distinct from $K'_{3,1}$). Also, $K'_{2,1}$ for which there are $6$ terms
is distinct from the other knots with $a>3$ in equation (\ref{eq:JKpb1}) and from the knots in equations
(\ref{eq:JKpaeqb}) and (\ref{eq:JKpaeqb2}) because the third term of $K'_{2,1}$ is $t^3$ whereas in these two equations
the power of the third term is even.

Consider the formulas with six terms i.e. equations (\ref{eq:JKpb1}), (\ref{eq:JKpaeqb}) and (\ref{eq:JKpaeqb2}).
In equation (\ref{eq:JKpb1}) the last term has coefficient $1$, which distinguishes it from the other two equations
with last coefficient $-1$.
For these two equations, suppose that $J_{K'_{b,b}}=J_{K'_{c+2,c}}$ for some integers $b>1$ and $c>1$.
By comparing the different terms in these polynomials we must have:

$2b=2c+2$ and $b^2+2b=c^2+4c$. From the first equation $b=c+1$ so $b^2+2b=c^2+4c+3\neq c^2+4c$.

Also, the knots $K'_{a,1}$ are distinguished by $a$, $K'_{b,b}$ are distinguished by $b$ and $K'_{b+2,b}$ are
also distinguished by $b$.

For the remaining subfamilies with eight terms the proof that they all have distinct polynomials $J$ is exactly like for $K_{a,b}$
in the previous case.

$\bullet$ Consider the family $K''_{a,b}$. For $J$ one has $4$, $5$, $6$, $7$ or $8$ terms.
The knots with $4$ terms are $K''_{2,1}$ (with second term $-t^3$), $K''_{3,1}$ (with second term $-t^4$)
and the knots $K''_{b+2,2}$ of equation (\ref{eq:JKppaeqb2}) (with second term $-t^{2b+2}$ and $b>1$ so $2b+2\ge 6$).
Thus, all these knots are distinct.

The only knot with $5$ terms is $K''_{4,1}$ so it is distinct from all other knots.

The knots with $6$ terms occur in equations (\ref{eq:JKppb1}) and (\ref{eq:JKppa2}). These two subfamilies are distinguished
by the coefficient of the last term ($-1$ and $1$). Also, $K''_{a,1}$ are distinguished by $a$ and $K''_{2,b}$ are distinguished by $b$.

The knots with $7$ terms occur in equations (\ref{eq:JKppaeqb1}) and (\ref{eq:JKppaeqb3}). Again, the two subfamilies are distinguished 
by the coefficient of the last term ($1$ and $-1$) and inside each family the knots are distinguished by $b$.

Finally, the knots with $8$ terms occur in equations (\ref{eq:JKppalb}) and (\ref{eq:JKppagb}). Once again, the two subfamilies are distinguished
by the coefficient of the last term ($1$ and $-1$). Inside each family the second term determines $a+b$ and the difference between the powers of the
seventh and the sixth terms determine $b-a$ or $a-b$. This together determines the couple $(a,b)$

For the second part of the theorem, it follows from Lemma~\ref{lemma_reducetoK} 
that if $h(K)=1$ then $K$ is a knot in one of the three listed families. One has to check that these knots do not have the Hopf crossing
number equal to $0$, i.e. that they are not right handed torus knots $T(n,n+1)$ (see Lemma~\ref{lemma_nocrossings}). 
From equation (\ref{eq:VT2}) the $J$ polynomials of $T(n,n+1)$ are equal to  $1-t^{n+1}-t^{n+2}+t^{2n+1}$. 
Consider all knots in the three families having $J$ polynomials with four terms. 
One checks easily that all these $J$ polynomials are different from the $J$ polynomials of the knots $T(n,n+1)$.
\end{proof}
\end{theorem}

Now we identify the knots $K$ with $h(K)=1$ and the classical crossing number $c(K)\le 10$ in the Rolfsen knot table\cite{R}.

\begin{proposition}
There are eight knots $K$ up to 10 crossings with $h(K)=1$, namely $K_{1,1}=4_1$,
$K_{2,1}=9_{42}$, $K'_{1,1}=3_1$,
$K'_{2,1}=5_2$, $K'_{3,1}=10_{132}$, $K'_{2,2}=10_{145}$,
$K''_{2,1}=\overline{5_1}=T(2,5)$ and $K''_{3,1}=10_{124}=T(3,5)$.
\begin{proof}
The span of the Jones polynomial is a lower bound on the number of crossings.
It is equal to the span of the polynomial $J$ minus $2$, as one divides by
$1-t^2$. Thus it is sufficient to consider the knots with span of $J$ not greater than $12$. 
Apart from the listed knots there are also the knots
$K_{3,1}$, $K_{2,2}$, $K'_{4,1}$, $K''_{4,1}$, $K''_{2,2}$ and $K''_{4,2}$.
Except for $K_{3,1}$, checking their Jones polynomial shows that these knots have more than 10
crossings. The knot $K_{3,1}$ was considered in example~\ref{ex:k31} in section~\ref{sec:examples}.
As it is distinct from $K_{2,1}=9_{42}$ and $\overline{9_{42}}$ and there are no knots up to 10 crossings sharing the Jones
polynomial with $9_{42}$, the knot $K_{3,1}$ has more than 10 crossings.

Using Lemma~\ref{lemma_boundcross} one checks that $K_{1,1}$, $K_{2,1}$, $K'_{1,1}$, $K'_{2,1}$ and $K''_{2,1}$
 have diagrams with at most $9$ crossings. Their Jones polynomials allow to identify these knots.
For $K'_{3,1}$ and $K'_{2,2}$  the lemma gives an upper bound of 13 crossings.
Drawing their diagrams one can easily reduce the number of crossings to 11. Using their DT codes and KnotFinder\cite{CL}
one identifies $K'_{3,1}$ as $10_{132}$ and $K'_{2,2}$ as $10_{145}$.

Finally $K''_{3,1}$ is the torus knot $T(3,5)=10_{124}$ (apply $\Om_\infty$ to the diagram of $T(3,5)$ in Figure~\ref{ex_Txx} 
to get $K''_{3,1}$).
\end{proof}
\end{proposition}

Let us study now the effect of taking the mirror image of a knot $K$ with $h(K)\le 1$. 
We show that, except for the knots $3_1$ and $4_1$, the Hopf crossing number of such a mirror image is greater than $1$.
The right handed trefoil knot satisfies $h(\overline{3_1})=h(T_{2,3})=0$ and $h(3_1)=h(K'_{1,1})=1$.
The figure eight knot $4_1=K_{1,1}$ is amphicheiral.

\begin{proposition}
Let $K$ be a knot with $h(K)\le 1$, $K\neq T(2,3)$, $K\neq K_{1,1}$ and $K\neq K'_{1,1}$. Let $\overline{K}$ be the mirror of $K$.
Then $h(\overline{K})>1$. In particular all such knots $K$ are chiral and none is a mirror image of another one. 
\begin{proof}
Notice that as $V_{\overline{K}}(t)=V_K(t^{-1})$, $J_{\overline{K}}(t)=-t^s J_K(t^{-1})$ where $s$ is the degree of $J_K(t)$.
Thus,  $J_{\overline{K}}$ is obtained from $J_K$ by symmetrically reversing the coefficients and multiplying by $-1$.
In particular the number of terms in $J_{\overline{K}}$ and $J_K$ is the same. Let us use in this proof the notation $g(K)$
for the ordered tuple of successive gaps between the powers of $J_K$. For instance, as the $J$ polynomials of $T(n,n+1)$ are equal to 
$1-t^{n+1}-t^{n+2}+t^{2n+1}$, we have $g(T(n,n+1))=(n+1,1,n-1)$. Obviously, $g(\overline{K})$ is obtained from $g(K)$
by reversing the order of the gaps.

Suppose that $K$ satisfies the hypothesis of the proposition. We need to show that $J_{\overline{K}}$ is not equal to any $J$ of
a knot with $h(K)\le 1$. We consider case by case $J_K$ with 4,5,...,8 terms.

$\bullet$ The knots $K$ with $J_K$ having 4 terms (excluding $K_{1,1}$ and $K'_{1,1}$) are $T(n,n+1)$ (with Hopf crossing number 0),
$K_{2,1}$, $K_{3,1}$, $K'_{3,1}$, $K''_{2,1}$ and $K''_{b+2,b}$, $b\ge 1$ (the exceptional case $K''_{3,1}$ can be included in this family). 
The gaps for $J$ of these knots are $g(K_{2,1})=g(K_{3,1})=(1,6,1)$, $g(K'_{3,1})=(1,3,3)$, $g(K''_{2,1})=(3,3,1)$ 
and $g(K''_{b+2,b})=(2b+2,2,2b)$. 

Reversing any of these gaps one never gets $(n+1,1,n-1)$ so these are not mirror images of 
torus knots $T(n,n+1)$. This is also true for the knots $T(n,n+1)$ as the equation $(n+1,1,n-1)=(n'-1,1,n'+1)$ has no solutions.
Thus $h(\overline{K})>0$ (notice that this is obviously also true for knots $K$ having $J_K$ with more than 4 terms).
Furthermore, this implies that $h(\overline{T(n,n+1)})>1$ for $n>2$.

From examples~\ref{ex:k21} and \ref{ex:k31} in section~\ref{sec:examples} it follows that $K_{2,1}$, $\overline{K_{2,1}}$, 
$K_{3,1}$ and $\overline{K_{3,1}}$ are pairwise distinct ($K_{2,1}=9_{42}$ is chiral). 
Thus, in particular, $h(\overline{K_{2,1}})>1$ and $h(\overline{K_{3,1}})>1$.
$K'_{3,1}$ is not the mirror image of $K''_{2,1}$ because the first knot is $10_{132}$ and the second is $\overline{5_1}$.
Finally, for $K''_{b+2,b}$ the last even $2b$ cannot be equal to $1$ or $3$ and the equation
$(2b+2,2,2b)=(2b',2,2b'+2)$ has no solutions. Thus $h(K)>1$ for all knots with $J_K$ having 4 terms.

$\bullet$ The knot $K''_{4,1}$ is the only knot with $J_{K''_{4,1}}$ having five terms and this polynomial is changed after reversing
the coefficients and multiplying by minus one. So $K''_{4,1}$ is chiral and $h(K''_{4,1})>1$.

$\bullet$ Consider the gaps of knots $K$ with $J_K$ having 6 terms. These are:
\[g(K'_{2,1})=(1,2,2,1,1)\]
\[g(K_{a,1})=(1,a+3,a-2,2,a-3),\; a\ge 4\]
\[g(K_{b,b})=(1,2b+1,1,b^2-2,2),\; g(K_{b+2,b})=(1,2b+3,1,b^2+2b-4,2),\; b>1\]
\[g(K'_{a,1})=(1,a,a,2,a-3),\; a\ge 4\]
\[g(K'_{b,b})=(1,2b-1,1,b^2-1,2),\; g(K'_{b+2,b})=(1,2b+1,1,b^2+2b-3,2),\;b>1\] 
\[g(K''_{a,1})=(a+1,2,a,a-4,1),\; a\ge 5\]
\[g(K''_{2,b})=(b+2,2,b-1,b,1),\;b>1\]
None of the gaps starts with $1,1$ so $h(\overline{K'_{2,1}})>1$.
We get the same inequality for $K_{b,b}$, $K_{b+2,b}$, $K'_{b,b}$ and $K'_{b+2,b}$ because the last gap for these knots is $2$
and none of the gaps starts with $2$.
The reversed gap of $K''_{a,1}$ is $(1,a-4,a,2,a+1)$ and the reversed gap of $K''_{2,b}$ is $(1,b,b-1,2,b+2)$.
None of these gaps can be equal to $(1,a'+3,a'-2,2,a'-3)$ or $(1,a',a',2,a'-3)$.
The last possibility of gaps symmetry to exclude is when $a=4$ so $g(K_{a,1})=(1,7,2,2,1)$ and $g(K'_{a,1})=(1,4,4,2,1)$ but
these gaps are not symmetric and reversing the first one does not give the second one.

$\bullet$ The two families of knots $K$ with $J_K$ having 7 terms have coefficient $2$ or $-2$ in the sixth term.
Thus $J_{\overline{K}}$ has $2$ or $-2$ in the second term so $h(\overline{K})>1$.

$\bullet$ Consider the gaps of knots $K$ with $J_K$ having 8 terms. These are:
\[g(K_{b+1,b})=(1,2b+2,1,b^2+b-3,1,1,1),\;b>1\]
\[g(K_{a,b})=(1,a+b+1,1,ab+b-a-2,a-b,2,a-b-2),\;b>1,a>b+2\]
\[g(K'_{b+1,b})=(1,2b,1,b^2+b-2,1,1,1),\;b>1\]
\[g(K'_{a,b})=(1,a+b-1,1,ab+b-a-1,a-b,2,a-b-2),\;b>1,a>b+2\]
\[g(K''_{a,b})=(a+b,2,a+b-2,ab-2b-1,1,b-a+1,1),\;a\le b,a>2\]
\[g(K''_{a,b})=(a+b,2,a+b-2,ab-a-b+1,1,a-b-3,1),\;b>1,a>b+3\]
None of the gaps starts with $1,1,1$ so $h(\overline{K_{b+1,b}})>1$ and $h(\overline{K'_{b+1,b}})>1$.
The reversed gaps of $K''_{a,b}$ end with $(...,a+b-2,2,a+b)$ (for both cases of $K''_{a,b}$) so
the last gap $a+b$ is larger than the gap $a+b-2$.
The gaps of $K_{a,b}$ and $K'_{a,b}$ end with $(...,a-b,2,a-b-2)$ so the last gap $a-b-2$ is smaller than
$a-b$. Thus $h(\overline{K''_{a,b}})>1$.
Finally the next to last gap of $K_{a,b}$ and $K'_{a,b}$ is 2 which cannot be the second gap for these knots.
\end{proof}
\end{proposition}

The Jones polynomials of knots $K$ with $h(K)=1$ have relatively simple form (the $J$ polynomials have at most 8 terms).
We show now that the coefficients of these Jones polynomials are bounded in absolute value by $2$.

\begin{proposition}
If $K$ is a knot with $h(K)=1$ then its Jones polynomial can have as coefficients only $\pm 1$ and $\pm 2$.
\begin{proof}
Consider the knots $K_{a,b}$.
If $a+b$ is even we can write the $J$ polynomials grouping terms in couples:
\[ J_{K_{a,b}}=(1-t^{a+b+2})-t(1-t^{a+b+2})+t^{ab+2b+1}(1-t^{a-b})-t^{ab+a+b+3}(1-t^{a-b-2}) \]
The Jones polynomial is obtained from $J$ by dividing by $1-t^2$ and multiplying by $t$ to some power
(the last operation only shifts the coefficients so we can ignore it).
Dividing the first couple $1-t^{a+b+2}$ by $1-t^2$ contributes $1$ to some coefficients of the Jones polynomial.
The same is true for the third couple, whereas for the second and the fourth dividing by $1-t^2$ contributes $-1$.
Adding these contributions we see that the coefficients are bounded in absolute value by $2$.

If $a+b$ is odd we use a different grouping:
\[ J_{K_{a,b}}=1-t^{ab+a+b+1}-t(1-t^{ab+2b})-t^{a+b+2}(1-t^{ab+a-b-1})+t^{a+b+3}(1-t^{ab}) \]
Again, we get $1$ from two couples and $-1$ from the other two couples, which gives the required bound
on the Jones polynomial coefficients.

For the two other families of knots the reasoning is similar so let us just group in couples the terms of $J$ polynomials.
For $K'_{a,b}$, suppose that $a+b$ is even:
\[J_{K'_{a,b}}=-1(1-t^{a+b})+t(1-t^{a+b})+t^{ab+2b}(1-t^{a-b})-t^{ab+a+b+2}(1-t^{a-b-2}) \]
Suppose now that $a+b$ is odd:
\[J_{K'_{a,b}}=-1(1-t^{ab+2b})+t(1-t^{ab+a+b-1})+t^{a+b}(1-t^{ab+2})-t^{a+b+1}(1-t^{ab+a-b-1}) \]
For $K''_{a,b}$, suppose that $a+b$ is even:
\[J_{K''_{a,b}}=1-t^{a+b}-t^{a+b+2}(1-t^{a+b-2})-t^{ab+a+b+1}(1-t^{a-b-2})+t^{ab+a+b+2}(1-t^{a-b-2}) \]
Suppose now that $a+b$ is odd:
\[J_{K''_{a,b}}=1-t^{ab+a+b+1}-t^{a+b}(1-t^{ab+2})-t^{a+b+2}(1-t^{ab+a-b-3})+t^{2a+2b}(1-t^{ab-2b}) \]
\end{proof}
\end{proposition}

\begin{proposition}\label{prop_torus}
Let $k>0$, then $h(T(n,n+k))=k-1$ for $k\le 3$ except for $h(T(2,5))=1$.
Also, $h(T(n,n+k))\ge 2$ for $k>3$.
\begin{proof}
For $k=1$ or $2$ the equality follows from Lemma~\ref{lemma_nocrossings} and section~\ref{sec:examples}.
From equation (\ref{eq:VT1}) the $J$ polynomial of $T(n,n+k)$ is $1-t^{n+1}-t^{n+k+1}+t^{2n+k}$.
If $k>2$ the $J$ polynomials are different from the $J$ polynomials of $T(n,n+1)$ and $T(n,n+2)$.
One checks case by case that they are also different from the $J$ polynomials with four terms of all other knots $K$
with $h(K)=1$ except for $K''_{2,1}$ which is $\overline{5_1}=T(2,5)$.
Thus, excluding $T(2,5)$, if $k>2$ then $h(T(n,n+k))\ge 2$ and, as in section~\ref{sec:examples} we saw that
$h(T(n,n+3))\le 2$, it follows that $h(T(n,n+3))=2$.
\end{proof} 
\end{proposition}

\section{Algebraic knots $K$ with $h(K)\le 1$}\label{sec:algknots}

An algebraic knot is an iterated torus knot of type $\{(p_1,q_1),(p_2,q_2),\ldots,(p_s,q_s)\}$.
Denote such a knot $K_s$.
It is constructed inductively: $K_1$ is the torus knot of type $(p_1,q_1)$, and, for $2\le i\le s$, $K_i$
is a $(p_i,q_i)$ cable of $K_{i-1}$.
Also, for every $i\in\{1,\ldots,s\}$, $gcd(p_i,q_i)=1$, $1<p_i<q_i$, and, for every $i\in\{1,\ldots,n-1\}$, $p_ip_{i+1}q_i<q_{i+1}$ (see \cite{Tg}).

In order to detect which knots $K$, with $h(K)\le 1$, are algebraic knots, we use a formula for the colored Jones polynomial of the $(p,q)$ cabling
of a knot $K$, denoted $K(p,q)$, from \cite{Tr} (see also \cite{Mo}):

\[V_{K(p,q)}(n)=A^{-pq(n^2-1)}\Sum_{j=-\frac{n-1}{2}}^\frac{n-1}{2}A^{4qj(jp+1)}V_K(2pj+1)\label{eq:Kpq}\tag{Kpq}\]

Note, that the conventions in \cite{Tr} are: $J$ for our $V$, $t$ for our $A$, $s$ for $p$ and $r$ for $q$.
The colored Jones $V_K(n)$ is normalized so that $V_U(n)=[n]$ for the unknot $U$, where:
\[[n]=\frac{A^{2n}-A^{-2n}}{A^2-A^{-2}}\]
The usual Jones polynomial of $K$, denoted $V_K$, equals $\frac{V_K(2)}{[2]}$ and $V_K(1)=1$ for every $K$.
Also, by convention, $V_K(-n)=-V_K(n)$.

We have,
\[[2](A^2-A^{-2})=A^4-A^{-4}=t^{-1}-t=t^{-1}(1-t^2)\]
Let $J_K(n)=(A^2-A^{-2})V_K(n)$. For $n=2$, it is the polynomial $J_K$ of section \ref{sec:classification}, up to some power of $t$.
Indeed, by definition of $J_K$, we have, for some $a\in\Z$:
\[J_K=t^a(1-t^2)V_K=t^a t [2](A^2-A^{-2})V_K=t^{a+1}J_K(2)\]

Let $K_s$ be an algebraic knot of type $\{(p_1,q_1),\ldots,(p_s,q_s)\}$.
We want to compute the highest power of $A$ in $J_{K_s}$ and a few gaps between the following decreasing powers of $A$.

Let $l_1=-2(q_1-1)(p_1-1)+2$ and, for $1<i\le s$, $l_i=-2q_ip_i+2q_i+p_il_{i-1}$.

\begin{lemma}
For any $i\ge 1$, $-2q_ip_i<l_i<0$
\begin{proof}
Clearly $-2q_1p_1<l_1<0$. Assume the lemma is true up to $i-1$.
\[l_i=-2q_ip_i+2q_i+p_il_{i-1}>-2q_ip_i+2q_i-2p_iq_{i-1}p_{i-1}>-2q_ip_i\]
The last inequality comes from $q_i>p_ip_{i-1}q_{i-1}$.
Also, $l_i$ is clearly negative by induction.
\end{proof}
\end{lemma}

\begin{lemma}\label{lemma_ls}
For $s\ge 2$, $4q_s-4(p_1+q_1)p_2\ldots p_s+2l_{s-1}>0$
\begin{proof}
First we check it for $s=2$:
\[4q_2-4(p_1+q_1)p_2-4(q_1-1)(p_1-1)+4>4q_1p_2p_1-4(p_1+q_1)p_2-4q_1p_1+4p_1+4q_1=\]
\[=4(p_2-1)q_1p_1-4(p_1+q_1)(p_2-1)=4(p_2-1)(q_1p_1-q_1-p_1)>0\]
Assume $s>2$. By the previous lemma and using $q_s>p_sp_{s-1}q_{s-1}$, it is sufficient to prove
\[4p_sp_{s-1}q_{s-1}-4(p_1+q_1)p_2\ldots p_s-4p_{s-1}q_{s-1}>0\]
or
\[(p_s-1)p_{s-1}q_{s-1}-(p_1+q_1)p_2\ldots p_s>0\]
As $\frac{p_s-1}{p_s}\ge\frac{1}{2}$, it is sufficient to show
\[\frac{(p_1+q_1)p_2\ldots p_{s-1}}{p_{s-1}q_{s-1}}<\frac{1}{2}\]
One checks easily by induction on $s>2$ that
\[\frac{(p_1+q_1)p_2\ldots p_{s-1}}{p_{s-1}q_{s-1}}\le\frac{p_1+q_1}{q_2} <\frac{2q_1}{p_1p_2q_1}\le \frac{1}{2}\]
\end{proof}
\end{lemma}

\begin{proposition}\label{prop_gaps}
Let $n\ge 1$. The highest power of $A$ in $J_{K_s}(n)$ is $(n-1)l_s+2$.
For $n\ge 2$, the three consecutive gaps between the four highest powers of $A$ are:

$4p_1\ldots p_s(n-1)+4$, $4(q_1-p_1)p_2\ldots p_s(n-1)$, $4p_1\ldots p_s(n-1)-4$.
\begin{proof}
If $n=1$, one has for any knot $K$, $J_K(1)=(A^2-A^{-2})V_K(1)=A^2-A^{-2}$ and the highest power of $A$ is indeed $2$.
For the rest of the proof we assume that $n\ge 2$.

Assume that $s=1$ so $K_1$ is a torus knot, $K_1=T(p,q)$, $p<q$ (to simplify the notation we use $p$ and $q$ instead of $p_1$ and $q_1$ in the case $s=1$).
Using $[n](A^2-A^{-2})=A^{2n}-A^{-2n}$ in formula (\ref{eq:Kpq}), with $K$ the unknot, we get:
\[J_{T(p,q)}(n)=A^{-pq(n^2-1)}\Sum_{j=-\frac{n-1}{2}}^\frac{n-1}{2}A^{4qj(jp+1)}(A^{4pj+2}-A^{-4pj-2})\]
For a given $j>0$, consider the four terms in the sum coming from $\pm j$:
\[A^{4qpj^2}(A^{4qj+4pj+2}-A^{4qj-4pj-2}-A^{-4qj+4pj-2}+A^{-4qj-4pj+2})\]
If $j\ge 1$, we check that the term with the lowest degree coming from $\pm j$ has higher degree than any of the terms coming from $\pm(j-1)$:
\[4qpj^2-4qj-4pj+2>4qp(j-1)^2+4q(j-1)+4p(j-1)+2\]
or
\[4(qj-1)(pj-1)-2>4(qj-q+1)(pj-p+1)-2\]
which is clearly true if $j\ge 1$.

Thus, the highest power of $A$ comes from $j=\frac{n-1}{2}$ and equals:
\[-pq(n^2-1)+pq(n-1)^2+2q(n-1)+2p(n-1)+2=(n-1)(-2pq+2q+2p)+2=\]
\[=(n-1)(-2(p-1)(q-1)+2)+2=(n-1)l_1+2\]
The first three gaps come from $j=\pm\frac{n-1}{2}$ and are:
$4p(n-1)+4$, $4(q-p)(n-1)$ and $4p(n-1)-4$.

Suppose now that $s>1$ and the proposition is true for all $K_i$, $i<s$.
Formula (\ref{eq:Kpq}) holds when we replace $V(n)$ with $J(n)$, since the second polynomial is $A^2-A^{-2}$ times the first.
Using induction on $s$ and this modified formula, the term for $j=\frac{n-1}{2}$ gives a power of $A$:
\[-p_sq_s(n^2-1)+p_sq_s(n-1)^2+2q_s(n-1)+((p_s(n-1)+1)-1)l_{s-1}+2=\]
\[=(n-1)(-2p_sq_s+2q_s)+p_s(n-1)l_{s-1}+2=(n-1)l_s+2\]
Also, for $j=\frac{n-1}{2}$, the highest gap in $J_{K_{s-1}}(2p_s j+1)=J_{K_{s-1}}(p_s(n-1)+1)$ is, by induction,
\[4p_1...p_{s-1}((p_s(n-1)+1)-1)+4=4p_1...p_{s-1}p_s(n-1)\]
Similarly, the two following gaps are: $4(q_1-p_1)p_2...p_s(n-1)$ and $4p_1...p_s(n-1)-4$.

We will prove now that the other terms in the sum in the formula for $J_{K_s}$ 
(i.e. for $j<\frac{n-1}{2}$) have sufficently low powers so they do not perturb the four highest powers coming from the term for $j=\frac{n-1}{2}$.
Denote by $d(j)$ the highest power of $A$ for the $j$ term in the sum (we ignore the common factor in front of the sum).

Let $j\ge\frac{1}{2}$. We compare $d(j)$ and $d(-j)$.
By induction the highest power of $A$ in $J_{K_{s-1}}(2p_sj+1)$ is $2p_sjl_{s-1}+2$. As
$J_{K_{s-1}}(2p_s(-j)+1)=-J_{K_{s-1}}(2p_sj-1)$, its highest power of $A$ is $(2p_sj-2)l_{s-1}+2$. The difference between these powers is $-2l_{s-1}$.
The coefficients in front of these two terms are $A^{4q_sj(jp_s+1)}$ and $A^{4q_s(-j)((-j)p_s+1)}$ giving a difference of $-8q_sj$.
Thus, the overall difference is
\[d(-j)-d(j)=-8q_sj-2l_{s-1}\]
This means that the drop in power of $A$ between $j$ and $-j$ is
\[d(j)-d(-j)=8q_sj+2l_{s-1}\ge 4q_s+2l_{s-1}>0\] 
the last inequality from Lemma~\ref{lemma_ls}. For $j=\frac{n-1}{2}$ the drop is $4q_s(n-1)+2l_{s-1}$.
We check that even after substracting the three first highest gaps for $j=\frac{n-1}{2}$, the drop is still positive. 
The sum of these three gaps is $4(p_1+q_1)p_2\ldots p_s(n-1)$. Using Lemma~\ref{lemma_ls} twice (once for $l_{s-1}<0$ together with $n-1\ge 1$),
\[4q_s(n-1)-4(p_1+q_1)p_2\ldots p_s(n-1)+2l_{s-1}\ge\]
\[(n-1)(4q_s-4(p_1+q_1)p_2\ldots p_s+2l_{s-1})>0\]
This shows that the term for $j=-\frac{n-1}{2}$ does not affect the four highest powers of $A$ coming from the term $j=\frac{n-1}{2}$.

Now let $j\ge 1$ be an index in the sum. Since, for $n=2$, $j=\pm\frac{1}{2}$, this means that $n\ge 3$. 
We compare $d(j)$ and $d(j-1)$. For $J_{K_{s-1}}$ the difference between the term for $j-1$ and $j$
is, by induction, 
\[((2p_s(j-1)+1-1)l_{s-1}+2)-((2p_sj+1-1)l_{s-1}+2)=-2l_{s-1}p_s\]
For the coefficients the difference is
\[4q_s(j-1)((j-1)p_s+1)-4q_sj(jp_s+1)=-4q_s(1+p_s(2j-1))\]
Thus the drop in power of $A$ is
\[d(j)-d(j-1)=4q_s(1+p_s(2j-1))+2l_{s-1}p_s>p_s(4q_s+2l_{s-1})>0\]
The last inequality comes from Lemma~\ref{lemma_ls}.

For $j=\frac{n-1}{2}$ the drop is $4q_s(1+p_s(n-2))+2l_{s-1}p_s$.
We check again that after substracting the three first gaps, the drop, say $D$, is still positive i.e.
\[D=4q_s(1+p_s(n-2))-4(p_1+q_1)p_2\ldots p_s(n-1)+2l_{s-1}p_s>0\]
We have:
\[\frac{D}{p_s(n-1)}=4q_s\frac{1+p_s(n-2)}{(n-1)p_s}-4(p_1+q_1)p_2\ldots p_{s-1}+\frac{2l_{s-1}}{n-1}\]

As $n\ge 3$, $\frac{n-2}{n-1}\ge\frac{1}{2}$, so $\frac{1+p_s(n-2)}{(n-1)p_s}\ge\frac{1}{2}$. Also $p_s\ge 2$, $\frac{2}{n-1}\le 1$ and $l_{s-1}<0$, so:
\[\frac{D}{p_s(n-1)}\ge \frac{1}{2}(4q_s-4(p_1+q_1)p_2\ldots p_{s-1}p_s+2l_{s-1})>0\]
the last inequality from Lemma~\ref{lemma_ls} again. Thus $D>0$.

In conclusion, we have shown that: the term for $j_0=\frac{n-1}{2}$ gives the wanted $4$ highest powers of $A$; the terms for $j_0-1$ and $-j_0$ do
not affect these $4$ powers; all other $j$ terms have lower powers than the term for $j_0-1$ (since $d(j)-d(j-1)>0$ if $j\ge 1$ and 
$d(j)-d(-j)>0$ if $j>0$). Thus, we get the wanted $4$ highest powers of $A$ for $J_{K_s}(n)$. 
\end{proof}
\end{proposition}

\begin{corollary}\label{cor_2torus}
If $K_s$ is an algebraic knot with $h(K_s)\le 1$, $K$ not a torus knot (i.e. $s\ge 2$), then $K_s$ is doubly iterated of type
$\{(p_1,p_1+1),(2,q_2)\}$ and it belongs to the family $K''_{a,b}$.
\begin{proof}
Recall that $t=A^{-4}$. Thus, the gaps between the highest powers of $A$ become gaps between the lowest powers of $t$ divided by $4$ for the $J$ polynomials.
As $J_{K_s}$ equals $J_{K_s}(2)$, up to a power of $t$, the first gap of $J_{K_s}$ equals the first gap of $J_{K_s}(2)$ which is $p_1\ldots p_s+1$.
We check the gaps for $K_{a,b}$, $K'_{a,b}$ and $K''_{a,b}$ in section~\ref{sec:classification}.   
As the first gaps of $J$ of $K_{a,b}$ and $K'_{a,b}$ are equal to $1$, these knots are not algebraic. 
The second gap of the $J_{K_s}$ is $(q_1-p_1)p_2\ldots p_s$. For $K''_{a,b}$, the second
gap equals $2$ (except for the torus knot $K''_{2,1}=T(2,5)$). 
Thus, if $K''_{a,b}$ is an algebraic knot, which is not a torus knot, then $s=2$, $q_1-p_1=1$ and $p_2=2$. 
\end{proof}
\end{corollary}

We want to find necessary conditions on $q_2$ in the preceding corollary in order to have $h(K_s)\le 1$.
Let $K$ be of type $\{(p,p+1),(2,q)\}$ and let $K'$ be the torus knot $T(p,p+1)$.
Using formula (\ref{eq:Kpq}),
\[V_K(2)=A^{-6q}(A^{4q}V_{K'}(3)+V_{K'}(-1))=A^{-6q}(A^{4q}V_{K'}(3)-1)\]
Thus,
\[J_K(2)=(A^2-A^{-2})V_K(2)=A^{-6q}(A^{4q}J_{K'}(3)-A^2+A^{-2})\]
In the proof of Proposition~\ref{prop_gaps} is a formula for $J_{T(p,q)}(n)$. Specializing to $K'$ and $n=3$,
\[J_{K'}(3)=A^{-8p(p+1)}\left(A^{4p(p+1)}(A^{8p+6}-A^2-A^{-6}+A^{-8p-2})+A^2-A^{-2}\right)\]
\[=A^{-4p^2+4p+6}(1-A^{-8p-4}-A^{-8p-12}+A^{-16p-8}+A^{-4p^2-12p-4}-A^{-4p^2-12p-8})\]
Thus,
\[J_K(2)=A^{-6q+4q-4p^2+4p+6}(1-A^{-8p-4}-A^{-8p-12}+A^{-16p-8}+\]
\[+A^{-4p^2-12p-4}-A^{-4p^2-12p-8}-A^{-4q+4p^2-4p-4}+A^{-4q+4p^2-4p-8})\]
Switching to $t$ and ignoring the factor in front, so that the lowest term is $1$, we get:
\[J_K=1-t^{2p+1}-t^{2p+3}+t^{4p+2}+t^{p^2+3p+1}-t^{p^2+3p+2}-t^{q-p^2+p+1}-t^{q-p^2+p+2}\label{eq:Jdbl}\tag{JK}\]
We get now a necessary condition on $q$ in order to have $h(K)\le 1$:

\begin{proposition}\label{prop_2torus}
Let $K$ be algebraic of type $\{(p,p+1),(2,q)\}$. If $h(K)\le 1$ then $q=2p(p+1)+1$ and $K=K''_{p+1,p}$.
\begin{proof}
Let $q_0=2p(p+1)+1$.
As $K$ is algebraic, $q\ge q_0$. We have,
\[q_0-p^2+p+1=p^2+3p+2\]
Thus, for $q_0$, there are $7$ increasing terms in equation (\ref{eq:Jdbl}) with the sixth term having coefficient $-2$.
If $q>q_0$, there are $8$ increasing terms in this equation.
We use the $J$ polynomials of $K''_{a,b}$ in section ~\ref{sec:classification} to determine which $K''_{a,b}$ may be equal to $K$.
Because of the five first terms of $J_K$, the only candidates are $K''_{b+1,b}$, $b\ge 2$ with $7$ terms in $J$ and
$K''_{a,b}$, $2<a\le b$, with $8$ terms. We want to reject the second possibility in order to reject $q>q_0$.

Notice that for $q_0$ the last term in equation (\ref{eq:Jdbl}) is $t^{p^2+3p+3}$ and if $q>q_0$ the power of the last
term is greater than $p^2+3p+3$.
If $K''_{a,b}=K$, then, comparing the second terms of $J$ in equations (\ref{eq:JKppalb}) and (\ref{eq:Jdbl}), we get $a+b=2p+1$.
As $a\le b$, one checks easily that $ab\le p(p+1)$. The last term of $J$ for $K''_{a,b}$ is $t^{ab+a+b+2}$ and we have,
\[ab+a+b+2\le p(p+1)+2p+1+2=p^2+3p+3\]
So $K''_{a,b}$ cannot be equal to $K$ if $q>q_0$.

Finally, for $J_K$, the first gap in equation (\ref{eq:Jdbl}) is $2p+1$, whereas for $J_{K''_{b+1,b}}$ it is $2b+1$.
Thus, if $K=K''_{b+1,b}$, then $p=b$, so $K=K''_{p+1,p}$.
\end{proof}
\end{proposition}

Notice: we cannot conclude from the last proposition, that the knots $K''_{p+1,p}$ are algebraic. It will be shown to be true in Theorem~\ref{th_fp}.

In the remaining part of this section we will show that, if $K$ is algebraic and $h(K)\le 1$, then $C_{alg}(K)=h(K)$. This
gives a positive answer, for knots with $h(K)\le 1$, to the problem of Fiedler (see the introduction). Because of Corollary~\ref{cor_2torus}
and Proposition~\ref{prop_2torus}, we have to show that if $K$ is a doubly iterated knot of type $\{(p,p+1),(2,2p(p+1)+1)\}$, then $C_{alg}(K)\le 1$.
We also need to check that for torus knots $K$ in Proposition~\ref{prop_torus} satisfying $h(K)\le 1$, we have $C_{alg}(K)\le 1$.

Let $\Se$ be a small sphere in $\C^2$, centered at the origin. Let $(x,y)\in \Se$, with $x=r_xe^{i\theta_x}$, $y=r_ye^{i\theta_y}$ and
$r_x^2+r_y^2=\epsilon^2$. Intersecting $\Se$ with complex lines through the origin gives the Hopf map $p:\Se\to \C P^1\approx S^2$, 
$p(x,y)=[x,y]$. 

We can visualize $p$ with the help of Figure~\ref{hopf_fibr}. $\Se$ is decomposed into an infinite family of tori and two $S^1$'s. 
For fixed $r_x$ and $r_y$ with $r_xr_y>0$, $(x,y)$ lies on a torus such as $\partial T$ (see Figure~\ref{hopf_fibr}). Let us choose the solid torus $T$ 
in such a way that $(x,y)\in\partial T$ if $r_x=r_y$ ($T$ is 'half' of $S^3$, the other half being $T'$).

If $\theta_x$ varies, $\theta_x\in[0,2\pi]$, we get the longitudes of the tori. Similarly, for $\theta_y\in[0,2\pi]$, we get the meridians of the tori. 
We choose the parametrization in such a way that, if $(x,y)\in D$ then $\theta_x=0$.
Increasing $r_x$ (and keeping $r_x^2+r_y^2=\epsilon^2$), we get increasingly thinner tori inside $T$.
The two $S^1$'s correspond to $(r_x,r_y)=(\epsilon,0)$ (the fiber going through the center of $D$) and 
$(r_x,r_y)=(0,\epsilon)$ (the fiber going through the center of $D'$). 

Now, $p(x,y)=[r_xe^{i\theta_x}, r_ye^{i\theta_y}]=[r_x,r_ye^{i(\theta_y-\theta_x)}]$.
If $r_x>0$ the map: $f: (x,y)\to (r_x,r_ye^{i(\theta_y-\theta_x)})$ is well defined. Let us restrict $f$ to $\{(x,y)\in\Se,r_x\ge r_y\}$.
Then $f$ maps the solid torus $T$ onto $D$ and, by definition, it represents the Hopf map $p$, restricted to $T$.
Thus, if $L$ is a link in $T$ and $f(L)$ generic in $D$, then an arrow diagram of $L$ is obtained from $f(L)$ by adding some arrows and 
information of under/over at the crossings. In particular, the number of crossings of this arrow diagram equals the number of crossings of $f(L)$.
Using $f:T\to D$ we can now estimate $C_{alg}(K)$ for algebraic knots $K$ with $h(K)\le 1$.

Let $X$ be a complex plane algebraic curve with a singularity at the origin of $\C^2$. Let $\epsilon<<1$ be small and let $K=X\cap\Se$.
Assume that, for some $p\in\N$, the Puiseux expansion of $X$ in $0$ (see~\cite{EN}) is:
\[y=x^\frac{p+1}{p}(1+x^\frac{1}{2p})\]
We use the notations of \cite{EN}, where $y=x^\frac{q_1}{p_1}(a_1+x^\frac{q_2}{p_1p_2})$, so $q_1=p+1,p_1=p,q_2=1$ and $p_2=2$.
It is shown there that $K$ is an iterated torus knot of type $\{(p_1,\alpha_1),(p_2,\alpha_2)\}$
with $\alpha_1=q_1$ and $\alpha_2=q_2+p_1p_2\alpha_1$ (Proposition 1A.1 in \cite{EN}).

As, in our case, $\alpha_2=1+2p(p+1)$, $K$ is of type $\{(p,p+1),(2,1+2p(p+1))\}$.
Notice, that $p=1$ is not excluded, in which case $K$ is a $(2,5)$ cable of the trivial torus knot of type $(1,2)$, i.e. $K$ is
the torus knot $T(2,5)$.

We can choose a parametrization of $X$ with $t=re^{i\theta}$:

$x=t^{2p}$

$y=t^{2p+2}(1+t)$

For $(x,y)\in K=X\cap\Se$, we have,
\[r_x^2+r_y^2=r^{4p}+r^{4p+4}|1+re^{i\theta}|^2=r^{4p}+r^{4p+4}(1+r^2+2r\cos\theta)=\epsilon^2\]

We assumed that $\epsilon$ is small. It can be checked easily that, for each $\theta\in [0,2\pi)$, there is a unique $r$ satisfying the preceding
equality. Denote it by $r_\theta$. It is also easily checked, that $r_\theta$ is increasing for $\theta\in[0,\pi]$ and decreasing for $\theta\in[\pi,2\pi]$.

Now $x=r_\theta^{2p}e^{2ip\theta}$ and $y=r_\theta^{2p+2}e^{i(2p+2)\theta}(1+r_\theta e^{i\theta})$. As $\epsilon<<1$ is arbitarily small, we can also assume that 
$r_\theta <<1$. Then it is clear that $r_x>r_y$, so $(x,y)\in T$ and:
\[f(x,y)=(r_\theta^{2p},r_\theta^{2p+2}e^{2i\theta}(1+r_\theta e^{i\theta}))\]

We see that for $\theta \in[0,2\pi]$, $r_\theta^{2p}$ will be increasing from $0$ to $\pi$, then decreasing from $\pi$ to $2\pi$.
Also, for $\theta \in[0,2\pi]$, the argument of the second coordinate of $f(x,y)$ goes from $0$ to $4\pi$, because of $e^{2i\theta}$, and it does
so monotonically. Indeed, $\frac{d}{d\theta}arg(e^{2i\theta})=2$, whereas it is easily checked that $|\frac{d}{d\theta}arg(1+r_\theta e^{i\theta})|<1$ for
$r_\theta<<1$. Thus, for $\theta\in[0,2\pi]$, $f(x,y)$ is a curve going twice around the center of $D$, first closer to it (as $r_x$ increases), then away from it
(as $r_x$ decreases). It follows that $f$ sends $K$ onto a curve with $1$ crossing. This means that $C_{alg}(K)\le 1$. 
From Lemma~\ref{lemma_nocrossings}, $h(K)\neq 0$. As $h(K)\le C_{alg}(K)$, we get $C_{alg}(K)=1$.

Finally, we consider the torus knots in Proposition~\ref{prop_torus} satisfying $h(K)\le 1$. These are the torus knots $T(n,n+1)$, $T(n,n+2)$ and
$T(2,5)$. The knot $T(2,5)$ has been considered above as a $(2,5)$ cable of the trivial knot, so $C_{alg}(T(2,5))=1$. 
For $T(n,n+k)$, $k=1$ or $2$, Fiedler states in \cite{F} (remark on page 260), that $C_{alg}(T(n,n+k))=k-1$, with $T(n,n+1)$ being realized
by the curve $X:z^n=w^{n+1}$ and $T(n,n+2)$ by a slight perturbation of the curve with equation $z^n=w^{n+2}$. The second case can be checked
using the following paramatrization for the perturbation of equation $z^n=w^{n+2}$: $z=t^n$ and $w=t^{n+2}(1+t)$.
With the same proof as for the doubly iterated torus knots above, one can show that this perturbation gives the knot $T(n,n+2)$ 
with a projection on the disk $D$ with $1$ crossing.

We resume the results of this section in the following theorem:

\begin{theorem}\label{th_fp}
The family of algebraic knots $K$ with $h(K)\le 1$ consists of: torus knots $T(n,n+1)$, if $h(K)=0$; torus knots $T(n,n+2)$, $T(2,5)$ or $K''_{p+1,p}$
if $h(K)=1$. The knots $K''_{p+1,p}$ are doubly iterated torus knots of type $\{(p,p+1),(2,2p(p+1)+1)\}$. Furthermore, for all such knots $C_{alg}(K)=h(K)$.
\begin{proof}
The only remaining thing to check is that the knots $K''_{p+1,p}$ are indeed algebraic. Let $K$ be an algebraic knot of type $\{(p,p+1),(2,2p(p+1)+1)\}$.
We have shown above that $C_{alg}(K)=1$. But $h(K)\le C_{alg}(K)$. By Proposition~\ref{prop_2torus}, $K$ must be equal to $K''_{p+1,p}$.
\end{proof}
\end{theorem}

We conclude with two remarks.

\begin{figure}\centering
\includegraphics{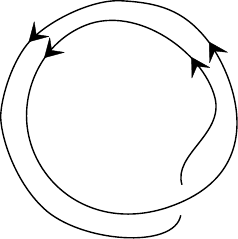}
\caption{Doubly iterated torus knot}
\label{ex_cable}
\end{figure}

\begin{remark}
It can be shown directly, using Reidemeister moves, that $K''_{p+1,p}$ is doubly iterated of type $\{(p,p+1),(2,2p(p+1)+1)\}$.
For example, we can transform $K''_{3,2}$ with $\Om_\infty$ and $\Om_4$ into a diagram shown in Figure~\ref{ex_cable}.
Then one could remove the two couples of arrows in this figure with appropriate Reidemeister moves to get a cable of the trefoil knot $T(2,3)$ and
check that it is the $(2,13)$ cable. This generalizes to any number $p$ of couples of arrows (for $p=1$ we get the middle diagram of Figure~\ref{ex_51}).
However, this proves only that such doubly iterated torus knots $K$ satisfy $h(K)=1$ and we have shown above that $C_{alg}(K)=1$, which is stronger.

It can also be checked that changing the crossing in Figure~\ref{ex_cable} gives a $(2,11)$ cable of the trefoil, which is also the knot $K''_{4,1}$.
More generally, we would get the family $K''_{p+2,p-1}$, $p\ge 2$, of doubly iterated knots of type $\{(p,p+1),(2,2p(p+1)-1)\}$. Notice,
that these are not algebraic knots.
\end{remark}

\begin{remark} Fiedler has shown in \cite{F}, Theorem 1, the following inequality for $K$ algebraic of type $\{(p_1,q_1),\ldots,(p_n,q_n)\}$, $q_1<2p_1$:
\[C_{alg}(K)\ge(q_1-p_1)p_2\ldots p_n-1\]
The right hand side of this inequality equals one if: $q_1-p_1=2$ and $n=1$ (torus knots $T(p_1,p_1+2)$) or $q_1-p_1=1$, $p_2=2$ and $n=2$ 
(type $\{(p_1,p_1+1),(2,q_2)\}$).
In the second case, it follows from Proposition~\ref{prop_2torus} that $C_{alg}(K)\ge h(K)>1$, if $q_2>2p_1(p_1+1)+1$. Thus, for such $q_2$, Fiedler's
inequality is not sharp. This can be compared to Theorem 2 in \cite{F}, which gives a sharper inequality then Theorem 1 for doubly iterated torus knots
with $p_1$ and $q_1$ both odd.
\end{remark}

\end{document}